\documentclass[12pt]{article}

\usepackage[utf8]{inputenc}
\usepackage[T1]{fontenc}
\usepackage[a4paper,margin=1in]{geometry}
\usepackage{amsmath,amssymb,amsthm,mathtools}
\usepackage{array,booktabs,longtable}
\usepackage{comment}
\usepackage{enumitem}
\setlist[enumerate]{itemsep=-2pt, topsep=2pt}
\usepackage{nicematrix}
\usepackage{float}
\usepackage[dvipsnames]{xcolor}
\usepackage{url}
\usepackage{csquotes}
\usepackage{listings}
\usepackage{hyperref}
\usepackage{cleveref}
\usepackage{caption}
\usepackage{authblk}

\usepackage[backend=biber,
style=trad-alpha,doi=false]{biblatex} %Imports biblatex package
\mathtoolsset{showonlyrefs,showmanualtags}
\numberwithin{equation}{section}

\hypersetup{
    colorlinks=true,
    linkcolor=orange!80!black,
    filecolor=magenta,
    urlcolor=RoyalBlue,
    citecolor=RoyalBlue
}

\usepackage{tcolorbox}
\newtcolorbox{mybox}{colback=blue!5!white,colframe=blue!75!black}

\theoremstyle{definition}
\newtheorem{definition}{Definition}[section]
\theoremstyle{theorem}

\newtheorem{proposition}[definition]{Proposition}
\newtheorem{corollary}[definition]{Corollary}
\newtheorem{theorem}[definition]{Theorem}
\newtheorem{conjecture}[definition]{Conjecture}

\newtheorem{lemma}[definition]{Lemma}
\newtheorem*{theorem*}{Theorem}

\newtheoremstyle{italicremark}
  {.5\topsep}{.5\topsep}{\normalfont}{}
  {\itshape}{.}{.5em}
  {\thmname{#1}\thmnumber{ \textit{#2}}\thmnote{ \normalfont(#3)}}

\theoremstyle{italicremark}
\newtheorem{remark}[definition]{Remark}

\theoremstyle{remark}
\newtheorem{examplex}[definition]{Example}
\newenvironment{example}
  {\pushQED{\qed}\examplex}
  {\popQED\endexamplex}

\newcommand{\R}{\mathbb R}
\newcommand{\Z}{\mathbb Z}
\newcommand{\C}{\mathbb C}

\let\P\relax
\newcommand{\P}{\mathbb P}
\newcommand{\PP}{\mathbb P}
\newcommand{\A}{\mathbb A}
\newcommand{\F}{\mathcal{F}}

\newcommand{\cO}{\mathcal O}
\newcommand{\cL}{\mathcal L}

\newcommand{\cS}{\mathcal S}
\newcommand{\cK}{\mathcal K}
\newcommand{\cQ}{\mathcal Q}
\newcommand{\cE}{\mathcal H}

\newcommand{\Pen}{\operatorname{Pen}}

\newcommand{\Gr}{\operatorname{Gr}}
\newcommand{\Sym}{\operatorname{Sym}}

\newcommand{\im}{\operatorname{im}}
\newcommand{\ev}{\operatorname{ev}}
\newcommand{\rk}{\operatorname{rank}}

\newcommand{\Set}[1]{\left\{#1\right\}}

\newcommand{\conv}{\operatorname{Conv}}
\newcommand{\wt}{\widetilde}
\newcommand{\ratto}{\dashrightarrow}

\newcommand{\pen}{\operatorname{Pen}}
\newcommand{\ogr}{\mathring{\operatorname{Gr}}}
\newcommand{\trop}{\operatorname{trop}}
\newcommand{\ms}{\mathcal{S}}

\newcommand{\joris}[1]{\textcolor{orange}{joris: #1}}
\newcommand{\Deltainv}{\Delta_{3,1}}

\newcommand{\vzero}{v_1}
\newcommand{\vone}{v_2}
\newcommand{\vtwo}{v_3}
\newcommand{\vthree}{v_4}
\newcommand{\aletter}{h}
\newcommand{\vh}{v_\aletter}
\newcommand{\anotherletter}{s}

\newcommand{\subwheelone}{U}
\newcommand{\subwheeltwo}{V}

\newcommand{\slopeJ}[1]{J_{#1}^s}
\newcommand{\slopeI}[1]{I_{#1}^s}

\DeclareMathOperator{\gr}{Gr}

\title{From MHV Gravity Amplitudes and Hodges map to Pentad Varieties}

\author[1]{Joris Koefler\textsuperscript{*}}
\author[2]{Emanuele Ventura}
\author[2]{Hongmiao Yu}
\affil[1]{\small{School of Natural Sciences, Institute for Advanced Study, 1 Einstein Drive, Princeton, USA}}
\affil[2]{\small{Politecnico di Torino, Dipartimento di Scienze Matematiche ``G.\ L.\ Lagrange'', Corso Duca degli Abruzzi 24, Torino, Italy}}

\date{}

\begin{document}

\maketitle

\begin{abstract}
    We study some geometric aspects behind the Hodges matrix appearing in the context of scattering amplitudes for the $\mathcal N=8$ supergravity.
    Through a vector bundle construction, we show that the Hodges matrix gives rise to a canonical rational map from the two-step flag variety to symmetric matrices, called the {\it Hodges map}.
    This construction allows to define the amplitude geometrically: its vanishing locus is a divisor -- the {\it MHV gravity divisor}.
    Inspired by the Hodges map, we introduce the notion of Hadamard ratio of two varieties.
    In the case relevant to gravity, this is the ratio of two Grassmannians, which we call the \emph{pentad variety}, connecting it to the pentad constraint from algebraic statistics.
    Moreover, we show that the pentad variety is closely related to the combinatorial algebraic geometry of  Martin's affine slope variety, and the  phirotopes of Below, Krummeck, and Richter-Gebert.
\end{abstract}


\begingroup
  \hypersetup{linkcolor=black}
  \tableofcontents
\endgroup

\section{Introduction}\label{sec:intro}

Tree-level scattering amplitudes in quantum field theories are rational functions on kinematic spaces encoding the probabilities of certain outcomes of particle interactions.
In recent years, remarkable insights have been gained by finding combinatorial and geometric structures in these amplitudes, in particular for $\operatorname{Tr}\varphi^3$ theory and $\mathcal{N}=4$ supersymmetric Yang--Mills theory; see \cite{lam2025combinatorialgeometryparticlephysics} for a survey.
As pointed out in~\cite{Gravituhedron} the geometric and combinatorial features of $\mathcal N=8$ supergravity~\cite{Supergravity} amplitudes, however, are still more elusive.
For example, these amplitudes will display double poles in soft limits, complicating the analysis.
In the first part of this work we present a geometric interpretation of these amplitudes.
This gives rise to the pentad varieties which we will study in the latter part from a purely algebro-geometric perspective.
Motivated by the notion of little group scalings for scattering amplitudes, we will reduce the study of the pentad variety to that of the affine slope variety as first defined in \cite{martin2006slopesdeterminednpoints}.

\medskip

Our first central object will be the reduced maximally helicity violating (MHV) tree-level gravity amplitude in this theory.
Remarkably, this amplitude has a compact determinantal formula originally due to Hodges~\cite{Hodges2012}.
Let $n\geq 5$, and let $\lambda_i,\widetilde\lambda_i\in\C^2$ be the columns of two rank-two spinor matrices satisfying momentum conservation $\sum_{i=1}^n\lambda_i\widetilde\lambda_i^T=0$.
Write $p_{ij}=\langle ij\rangle\coloneqq\det(\lambda_i,\lambda_j)$ and $q_{ij}=[ij]\coloneqq\det(\widetilde\lambda_i,\widetilde\lambda_j)$, for the Pl\"ucker coordinates on $\gr(2,n)$.
Assume that $p_{ij}\neq0$ for all $i<j$.
The \emph{Hodges matrix} is the symmetric $n\times n$ matrix $\Phi$ with entries
\[
    \Phi_{ij}=
    \begin{cases}
        \displaystyle\frac{q_{ij}}{p_{ij}}, & i\neq j,\\[6pt]
        \displaystyle-\sum_{k\neq i}\frac{q_{ik}}{p_{ik}}
        \frac{\langle kx\rangle\langle ky\rangle}
        {\langle ix\rangle\langle i y\rangle}, & i=j,
    \end{cases}
\]
where $x,y\in\C^2$ are reference spinors chosen so that the denominators are nonzero.
%Momentum conservation makes $\Phi$ independent of this choice, and its generic rank is $n-3$.
Hodges' formula gives the reduced MHV gravity amplitude then as
\[
    \widetilde M_n\coloneqq \pm \frac{\det\Phi_{ijk}^{rst}}{p_{rs}p_{rt}p_{st}p_{ij}p_{ik}p_{jk}},
\]
where $i<j<k$ and $r<s<t$, and $\Phi_{ijk}^{rst}$ is obtained by deleting rows $\{i,j,k\}$ and columns labeled by $\{r,s,t\}$ (the correct sign in front depends on the chosen indices of rows and columns); see also~\cite{Gravituhedron}.
By the results in~\cite{SH}, pairs of spinor brackets $(p_{ij},q_{ij})$ satisfying momentum conservation and the Pl\"ucker relations, form the spinor-helicity variety $\mathrm{SH}(2,n,0)$, which is in turn isomorphic to the two-step flag variety $\F_n\coloneqq\operatorname{Fl}(2,n-2;n)$.
Hence $\widetilde{M}_n$ is a rational expression on $\F_n$.
It has been a long-standing open issue to find the geometric origin of Hodges' formula.
Our first main result is to give such a geometric interpretation by an intrinsic vector bundle construction on $\F_n$.
\begin{remark}
The reduced MHV gravity amplitude can also be written as 
\[
\widetilde{M}_n = \frac{N_n}{\prod_{i<j}p_{ij}},
\]
where a priori $N_n$ is a rational function on $\mathrm{SH}(2,n,0)$. Upon defining 
on the coordinate ring $Q_n$ of $\mathrm{SH}(2,n,0)$ a physically-motivated $\Z^{n+2}$-valued multidegree, in \cite{Gravity} the authors prove that $N_n$ is a polynomial class in $Q_n$ of a specific multidegree $d(n)\in \Z^{n+2}$. They conjectured that $N_n$ is the unique (up to scaling) polynomial of multidegree $d(n)$ being in the ideal of $Q_n$ given by the intersections of all monomial ideals $(p_{ij},q_{ij})\subset Q_n$. We recently disproved this conjecture for every $n\geq 7$ \cite{KoeflerVenturaYu2026}. 
\end{remark}

\medskip

We will show that the local expressions for Hodges matrix $\Phi$ glue to a rational map from $\F_n$ to the projectivized space of symmetric $n\times n$ matrices $\P(\Sym_\C(n))$:
\[
    H:\F_n\dashrightarrow\P(\Sym_{\C}(n)),\quad (L,W)\longmapsto[\Phi(L,W)],
\]
which we will call the \emph{Hodges map}.
Since it is also well known that Hodges matrix is generically of rank $n-3$, see~\cite{Hodges2012}, Hodges map $H$ induces another rational map $\widetilde H:\F_n\dashrightarrow \Sigma\subset \P(\Sym_{\C}(n))\times\gr(3,n)$ by sending $(L,W)\mapsto([\Phi(L,W)],\mathcal K_{(L,W)})$, where $\mathcal{K}_{(L,W)}$ is the kernel of $\Phi(L,W)$, and $\Sigma$ is the incidence variety defined as
\[
    \Sigma\coloneqq\left\{([A],K)\in\P(\Sym_{\C}(n))\times\gr(3,n)\mid K\subseteq\ker A\right\}.
\]
This incidence $\Sigma$ is isomorphic to the blow-up of the variety of symmetric $n$ matrices of rank at most $n-3$ in its singular locus, see Corollary~\ref{cor: blowup}, whose exceptional locus we denote by $E$.
Then, our main result gives a geometric interpretation of Hodges formula.

\begin{theorem}[The MHV gravity divisor] \label{thm:main_gravity}
    Let $n \geq 5$ and $\F_n^\circ$ be the Hodges locus in $\F_n$ where all Pl\"ucker coordinates $p_{ij}$ do not vanish.
    \begin{enumerate}
        \item The \emph{(open) gravity divisor} is the effective Cartier divisor $D_H^{\circ}\coloneqq \wt H^*E$ on the Hodges locus $\F_n^\circ$.
        \item For choices of $\lbrace i<j<k\rbrace$ and $\lbrace r<s<t\rbrace$, locally the open gravity divisor $D_H^\circ$ coincides with the vanishing locus in $\F_n^\circ$ of Hodges formula
        \[
            \Delta_H\coloneqq \pm \frac{\det\Phi^{rst}_{ijk}}{p_{ij}p_{ik}p_{jk}\,p_{rs}p_{rt}p_{st}}.
        \]
    \end{enumerate}
\end{theorem}

\medskip

The construction of the Hodges matrix with off-diagonal entries $q_{ij}/p_{ij}$ and of the Hodges map leads us to a second geometric object, which is an instance of what we will call a \emph{Hadamard ratio}.
Analogously to Hadamard products $X*Y$ of two projective varieties, \emph{Hadamard ratios} $X:Y$ are defined by taking coordinatewise ratios of $X$ and $Y$.
In the context of MHV gravity amplitudes, we focus on the Hadamard ratio of two Grassmannians of spinors, which is the Zariski closure of the image of
\[
    \Psi_n:
    \mathrm{Gr}(2,n)\times \mathrm{Gr} (2,n)\dashrightarrow
    \mathbb P^N,\quad (q_{ij},p_{ij})_{i<j}\longmapsto \left[\frac{q_{ij}}{p_{ij}}\right]_{i<j},
\]
where $N= \binom{n}{2}-1$, and the ratios $q_{ij}/p_{ij}$ are precisely the off-diagonal entries of Hodges matrix $\Phi$.
This notion generalizes in particular the idea of reciprocal spaces, which have attracted attention in recent years; see e.g.~\cite{Kummer_2019}.
For $n=5$, with coordinates on $\P^9$ denoted by $z_{ij}$, $\gr(2,5):\gr(2,5)$ is a $8$-dimensional hypersurface cut out by the classical \emph{pentad}
\[
\begin{matrix}
         \text{pen}_{12345}=&\quad z_{12} z_{13} z_{24} z_{35} z_{45}
        -z_{12} z_{13} z_{25} z_{34} z_{45}
        -z_{12} z_{14} z_{23} z_{35} z_{45}
        +z_{12} z_{14} z_{25} z_{34} z_{35} \\
        &+\,z_{12} z_{15} z_{23} z_{34} z_{45}
        -z_{12} z_{15} z_{24} z_{34} z_{35}
        +z_{13} z_{14} z_{23} z_{25} z_{45}
        -z_{13} z_{14} z_{24} z_{25} z_{35} \\
        &\,-\,z_{13} z_{15} z_{23} z_{24} z_{45}
        +z_{13} z_{15} z_{24} z_{25} z_{34}
        +z_{14} z_{15} z_{23} z_{24} z_{35}
        -z_{14} z_{15} z_{23} z_{25} z_{34},
    \end{matrix}
\]
appearing in algebraic statistics \cite{pentad}.
This coincidence suggests the Hadamard ratio $\gr(2,n):\gr(2,n)$ should be viewed as a higher pentad variety; hence, we refer to it as the \emph{pentad variety} $\Pen_n\coloneqq \gr(2,n):\gr(2,n) \subset \P^N$.
Moreover, we will see in Proposition \ref{prop: lifting to rank 2 symm matrices n=5} that, for generic data and up to a normalization factor, the pullback of any $2\times 2$ minor of the symmetric $5\times5$ matrix $Z_5$ with off-diagonal entries $q_{ij}/p_{ij}$ and diagonal entries filled in such that $\rk Z_5\leq 2$ returns the MHV gravity amplitude $\widetilde M_5$.
This connects the study of pentad varieties with that of these amplitudes.

\medskip
The remainder of our work is to study the pentad variety from an algebro-geometric point of view.
We investigate its dimension, defining equations, and its connection to further combinatorial algebraic varieties.
%We will also use tropical methods \cite{tropGeom}.
Let $J^p_n$ be the ideal in the polynomial ring $\C[z_{ij}\mid 1 \leq i<j\leq n]$ generated by all pentad constraints $\operatorname{pen}_\sigma$ on five-element subsets $\sigma$ of $[n]\coloneqq \{1,\ldots, n\}$.
Moreover, let $I^p_n$ be defined as the saturation of $J^p_n$ with respect to the product $\prod_{i<j}z_{ij}$ of the $N+1$ homogeneous coordinates on $\P^N$.
Then, our first main result about $\Pen_n$ can be summarized as follows.

\begin{theorem}\label{thm:main_pentad}
    The pentad variety $\pen_n\subset \P^N$ is irreducible of dimension $3n-7$ with prime vanishing ideal $\sqrt{I^p_n}$.
    In particular,
    \[
        \mathbb{V}\left(\sqrt{I^p_n}\right) = \mathbb{V}(I^p_n) = \gr(2,n):\gr(2,n)=\pen_n.
    \]
    Moreover, for $n\leq 10$, the ideal $I^p_n$ is radical, hence $I^p_n$ is the prime vanishing ideal. 
\end{theorem}

\medskip

Our third major object arises from the notion of \emph{little-group scaling} invariant variables, which are natural for the scattering amplitudes.
In our setup this little group scaling takes the following form:
Write the columns of generic representatives of the two Grassmannian $\gr(2,n)$ as $\nu_i(1,x_i)^T$ and $\mu_i(1,y_i)^T$ with $\nu_i,\mu_i\in \C^\times$, respectively.
Their Hadamard ratios are
\[
    z_{ij}=\frac{q_{ij}}{p_{ij}}
    =\frac{\mu_i\mu_j}{\nu_i\nu_j}\frac{y_j-y_i}{x_j-x_i}
    =t_it_j\frac{y_j-y_i}{x_j-x_i},
    \quad \text{where } t_i\coloneqq\frac{\mu_i}{\nu_i}.
\]
The little group scalings $t_it_j$ parametrize the dense torus of the toric variety $X_{\Delta}\subset \P^N$ associated to the second hypersimplex $\Delta\coloneqq \conv\{e_i+e_j\mid 1\leq i<j\leq n\}$.
The factors $(y_j-y_i)/(x_j-x_i)$ are slopes determined by points $(x_i,y_i)\in\C^2$.
This connects the pentad variety to the \emph{affine slope variety $\ms_n$}, as first defined in~\cite{martin2006slopesdeterminednpoints}, given as the Zariski closure in $\C^{N+1}$ with $N=\binom{n}{2}-1$ of the slope vectors of $n$ points in the complex plane with $x_i\neq x_j$ for $i\neq j$.
The following result makes this connection explicit.
% Here, take $n$ points $(x_i,y_i)\in \C^2$ such that none of the slopes $m_{ij}=\frac{y_j-y_i}{x_j-x_i}$ are infinite.
% Then, the affine slope variety is the Zariski closure of the slope vectors
% \[
%     \ms_n = \overline{\left\{(m_{12},\ldots, m_{n-1,n})\:\bigg\vert\: m_{ij}=\frac{y_j-y_i}{x_j-x_i} \text{ and } x_j\neq x_i\right\}}\subset \C^{\binom{n}{2}}.
% \]

\begin{corollary}\label{cor:main_aff_to_pen}
    The pentad variety $\pen_n$ is given by the Hadamard product of the projectivization of the affine slope variety $\P(\ms_n)$ with $X_{\Delta}$.
    Put in symbols, 
    \[
        \pen_n = \P(\ms_n)*X_{\Delta},
    \]
    as subvarieties of $\P^N$.
\end{corollary}

\begin{remark}
    The affine slope variety can also be defined by polynomials associated to graphs, see \cite{martin2006slopesdeterminednpoints}, thereby reducing the study of the pentad variety to a problem in combinatorial commutative algebra.
    This allows us in Section \ref{sec:affine_slope} to use these polynomials to study the defining equations of pentad varieties.
\end{remark}

Lastly, we show that the slope variety has an interpretation in terms of \emph{phirotopes}, the complex phase analogues of chirotopes introduced in~\cite{phirotopes}.
A \emph{realizable} rank-two phirotope records the phases of the $2\times2$ determinants of a complex vector configuration $Z_0,\ldots, Z_n\in \C^2$.
For a \emph{uniform} phirotope these phases are nonzero on every pair of distinct labels $i\neq j$.
For uniform realizable phirotopes, these phases can be interpreted, up to normalization, as slope vectors of affine points, hence connecting the affine slope variety to phirotopes.
The following result makes this precise.

\begin{proposition}
    For a uniform rank $2$ phirotope $\varphi$ with index set $E=\{0,1,\ldots,n\}$, denote its normalized phases by $\widehat{\varphi}(i,j)\coloneqq \varphi(i,j)/(\varphi(0,i)\varphi(0,j))$ for all $1\leq i<j\leq n$.
    Let $\widehat{\varphi}^{\,2}\coloneqq (\widehat{\varphi}(1,2)^2,\ldots,\widehat{\varphi}(n{-}1,n)^2)\in \C^{\binom{n}{2}}$ be the vector of their squares.
    Define $\mathcal{R}_n\subset \C^{\binom{n}{2}}$ as the set of all such points obtained from uniform realizable rank $2$ phirotopes $\varphi$ on $E$, that is
    \[
        \mathcal{R}_n \coloneqq \left\{ \widehat{\varphi}^{\,2} \mid \varphi \text{ is a uniform realizable rank $2$ phirotope on } E \right\} \subset \C^{\binom{n}{2}}.
    \]
    Then, the complex Zariski closure of $\mathcal{R}_n$ coincides with the affine slope variety $\ms_{n}$.
\end{proposition}

\noindent\textbf{Structure of the paper.}
Section \ref{sec:preliminaries} recalls the spinor-helicity formalism, and the identification of the spinor-helicity variety with the two-step flag variety $\mathcal F_n$.
Moreover, we introduce the notion of Hadamard ratios.
In Section \ref{sec:hodges-gravity-divisor}, we put the Hodges matrix in the context of vector bundles. We start by constructing the Hodges bundle, and a generalized Hodges matrix, define the Hodges map, and identify the associated vanishing determinant with the MHV gravity divisor, leading to  Theorem \ref{thm:main_gravity}.
Section \ref{sec:pentads} studies the Hadamard ratio map and the resulting pentad varieties from an algebro-geometric perspective. The proof of the main result will use the connection of pentad varieties with the affine slope varieties (see Theorem \ref{thm:main_pentad} and Corollary \ref{cor:main_aff_to_pen}). 
We end this section by showcasing another surprising incarnation of the affine slope variety as the locus of phases of realizable (uniform) phirotopes.

\medskip

\noindent \textbf{Acknowledgments.}
We are grateful for several discussions with Bernd Sturmfels, who also originally pointed out the observation that the pentad constraint gives the Hadamard ratio $\gr(2,5):\gr(2,5)$, as well as the potential connection to the phirotopes.
Moreover, we are very grateful to Ben Hollering, whose help and insights into the computational aspects of this work have been invaluable.
We thank Hannah Friedman for pointing out the connection formalized in Proposition \ref{prop:toric_times_slope_pentad}. 
Lastly, we thank Dario Antolini for a helpful discussion on Hadamard products.
J. K. is grateful for the support received from the Institute for Advanced Study and an anonymous donor. H.Y. wishes to express sincere gratitude to the Max Planck Institute for Mathematics in the Sciences (MPI MIS) for its hospitality during the preparation of this work. E.V. and H.Y. are members of the GNSAGA group of INdAM. E.V. and H.Y. acknowledge that the preparation of this article was partially funded by the Italian Ministry of University and Research, through the project {\it Applied Algebraic Geometry of Tensors}, PRIN 2022 Protocol no. 2022NBN7TL.

\medskip

\noindent \textbf{AI Disclosure.}
AI was used as a search engine and for proofreading purposes of drafts of this work.
Moreover, it was used for the implication (iii)$\Rightarrow$(i) in Corollary \ref{cor:realizability}, and for the pairing in (iii) of Lemma \ref{lem: induced det}.
The text of the article is entirely written by the authors, and they take full responsibility for its content, including all mathematical results. 

\medskip
\section{Preliminaries} \label{sec:preliminaries} 

In this section we introduce the spinor formalism, which is the language most frequently used in the study of MHV gravity scattering amplitudes.

\subsection*{Spinor Formalism}

We start by fixing two copies of the Grassmannian $\Gr(2,n)$ with Pl\"ucker coordinates denoted by $\langle ij\rangle$ and $[ij]$ where $1\leq i<j\leq n$, respectively.
So in matrix coordinates, these are $2\times 2$ minors of the corresponding $2\times n$ matrices representing points in the Grassmannian. 
Let $R_n \coloneqq \mathbb{C}[[ij],\langle ij\rangle \mid 1\leq i<j\leq n]$ be the polynomial ring generated by these Pl\"ucker coordinates. The variables $[ij]$ and $\langle ij\rangle$ are called {\it spinor variables}.
Let $p_i\in \C^4$ be the $4$-complex momentum of a massless particle $i$.
We say a particle is massless if $p_i\cdot p_i = p_i^2=0$, where the product is the Lorentzian inner product on $\mathbb{C}^4$.
One can define an isomorphism between $\mathbb{C}^4$ and $2\times 2$ complex matrices $p_i\mapsto P_i$ such that the vanishing of the inner product $p_i^2=0$ is equivalent to the condition $\det(P_i)=0$.
Since $P_i$ is a rank-one matrix,  one can write $P_i = \lambda_i\widetilde{\lambda_i}$, for column and row vectors $\lambda_i,\widetilde{\lambda_i}\in \mathbb{C}^2$.
The spinor variable $[ij]$ is the determinant of the $2\times 2$ matrix whose rows are $\widetilde{\lambda}_i,\widetilde{\lambda}_j$; the spinor variable $\langle ij\rangle$ is the determinant of the $2\times 2$ matrix whose columns are $\lambda_i, \lambda_j$. 

% \subsection*{Momentum Conservation}

% Next, we introduce momentum conservation i.e. $\sum_{i=1}^n P_i = 0$, where the equality is as $2\times 2$-matrices, in the language of spinors.
% Since $P_i$ is a rank-one matrix, as before we may write $(P_{i})_{\alpha\dot{\alpha}} = \lambda_{i\alpha}\tilde{\lambda}_{i\dot{\alpha}}$, using the index notation. 
% For any pair of indices $a,b\in [n]$, using the usual notation for contractions, we have 
% \[ 
%     \langle a i \rangle = \lambda_a^{\alpha}\lambda_{i\alpha} \mbox{ and } 
%     [i b] = \lambda_i^{\dot{\alpha}}\lambda_{b\dot{\alpha}}.
% \]
% Then 
% \[
%     \langle a i\rangle[i b] = \lambda_a^{\alpha}(\lambda_{i\alpha}\widetilde{\lambda}_i^{\dot{\alpha}})\widetilde{\lambda}_{b\dot{\alpha}} = \lambda_a^{\alpha}P_{i\alpha}^{\dot{\alpha}}\widetilde{\lambda}_{b\dot{\alpha}}.
% \]
% The condition $\sum_{i=1}^n P_i=0$ implies that $\sum_{i=1}^n P_{i\alpha}^{\dot{\alpha}} = 0$. Summing over $i$ in the expression above, we have
% \[
%     \sum_{i=1}^n \langle a i \rangle[i b] = \lambda_a^{\alpha}\left(\sum_{i=1}^n P_{i\alpha}^{\dot{\alpha}}\right)\widetilde{\lambda}_{b\dot{\alpha}} = 0. 
% \]
% These are $n^2$ (for each choice of $a,b\in [n]$) momentum conservation relations. 

\subsection*{Spinor-Helicity Variety}

From now on, we let $p_{ij}=\langle ij\rangle$ and  $q_{ij}=[ij]$ for $1\leq i<j\leq n$ and denote the polynomial ring in these variables also by $R_n =\C[p_{ij},q_{ij} \mid 1\leq i<j\leq n]$.
We define the ideal $L_n\subset R_n$ generated by the Pl\"ucker relations
\begin{align*}
        p_{ij}p_{k\ell}-p_{ik}p_{j\ell}+p_{i\ell}p_{jk},\\
        q_{ij}q_{k\ell}-q_{ik}q_{j\ell}+q_{i\ell}q_{jk},
\end{align*}
for all $1\leq i<j<k<\ell\leq n$, together with the momentum-conservation equations
\[
M_{ab}\coloneqq \sum_{r=1}^{n}p_{ar}q_{rb}=0,
        \qquad 1\leq a,b\leq n.
\]
The ring $Q_n=R_n/L_n$ is the coordinate ring of the spinor-helicity variety $\mathrm{SH}(2,n,0)$
of El Maazouz-Pfister-Sturmfels \cite{SH}. 

\medskip

\subsection*{Two-step Flag Varieties}

Lastly, we want to explicitly connect this construction with a flag variety.
Let $V=\C^n$ with its standard basis $e_1,\ldots,e_n$ and its standard nondegenerate symmetric bilinear form.  We look at the product of two Grassmannians $\Gr(2,n)\times \Gr(2,n)$. A point of the first Grassmannian is a $2$-dimensional plane $L\subset V$, 
with Pl\"ucker coordinates $p_{ij}$. A point of the second factor Grassmannian is another $2$-dimensional plane $M\subset V$, 
with Pl\"ucker coordinates $q_{ij}$.

Choose two $2\times n$ matrices $X$ and $Y$ whose row spans are $L$ and $M$. 
Then the momentum-conservation equations are the vanishing matrix multiplication $XY^T=0$. 
Equivalently, every vector in $L$ is orthogonal to every vector in $M$. 
Letting $W=M^\perp$, the condition $L\perp M$ becomes $L\subset W$. Hence a spinor-helicity point is the same datum as a two-step flag $L\subset W\subset V, \dim L=2, \dim W=n-2$. 
One can show that this identification yields an isomorphism between the spinor-helicity and the two-step flag variety $\mathcal F_n = \mathrm{Fl}(2,n-2;n)$ of partial flags $L\subset W\subset V$ \cite{SH}.

\subsection*{Hadamard Ratios}

In this section we introduce Hadamard ratios, which are a natural extensions of Hadamard products.
Let $X$ and $Y$ be two irreducible subvarieties of $\P^n$ with coordinates labeled by $x_i$ and $y_i$ respectively.
Then, the \emph{Hadamard product} $X*Y$ is defined as the Zariski closure of the image of the map
\[
    X\times Y \dashrightarrow \P^n,\quad  (x_i,y_i)\mapsto  x_iy_i.
\]
Hadamard products are well-studied objects, see for example \cite{HadamardProducts}.
In this work, we take a slightly different route.

\begin{definition}\label{def:Hadamard_ratio}
    Let $X$ and $Y$ be two irreducible subvarieties of $\P^n$, with coordinates labeled by $x_i$ and $y_i$, respectively.
    Then, we define the \emph{Hadamard ratio map} by
    \begin{align}
        \Psi:X\times Y\dashrightarrow \P^n,\quad (x_i,y_i)\mapsto  x_i/y_i,
    \end{align}
    and refer to the Zariski closure of the image $\overline{\im\Psi}$ as the \emph{Hadamard ratio} $X:Y$ between $X$ and $Y$.
\end{definition}

\begin{example}
    Suppose $Y=L\subset \P^n$ is a $d$-dimensional linear subspace.
    Then, for the variety $p=[1:1:\cdots:1]\in \P^n$, the Hadamard ratio $p:L$ is precisely the reciprocal linear space $L^{-1}$, which has attracted much attention in the literature; see e.g. \cite{Kummer_2019}. 
\end{example}

Denote by $T^n\subset \P^n$ the \emph{algebraic torus}, which is the set of all points in $\P^n$ where none of the coordinates vanishes.
Then, 
%if $X$ and $Y$ are subvarieties of $T^n$, then we get the following result, which is
we get natural extensions of a result in \cite[Section 5]{HadamardProducts}.

\begin{lemma}\label{lem:tropicalDiff}
    Let $X$ and $Y$ be two irreducible complex subvarieties of $\P^n$ and suppose that $X^{\circ}\coloneqq X\cap T^n\neq \emptyset$ and $Y^{\circ}\coloneqq Y\cap T^n\neq \emptyset$. 
    Then tropicalization of the Hadamard ratio $X : Y$ is the Minkowski difference of the tropicalizations of $X$ and $Y$. Put in symbols
    \[
        \trop((X:Y)) = \trop X - \trop Y.
    \]
\end{lemma}
\begin{proof}
    The Hadamard ratio map $\Psi$ is a monomial map, hence by \cite[Corollary 3.2.13]{tropGeom} the tropicalization of the image of $\Psi$ is given by image of $\trop (X^{\circ} \times Y^{\circ})$ under the tropicalization of $\Psi$.
    The latter is given by coordinate-wise subtraction,
    hence the statement.
\end{proof}

Hence, \cite[Proposition 5.4]{HadamardProducts} directly generalizes, which we record here in a slightly weaker form.
\begin{lemma}\label{lem:tropDim}
    Let $X$ and $Y$ be two irreducible subvarieties in $\P^n$.
    Denote by $H \subset (\mathbb{C}^*)^{n+1}/\mathbb{C}^*$ the largest subtorus acting on both
    $X$ and $Y$.
    %and $G \subset (\mathbb{C}^*)^{n+1}/\mathbb{C}^*$ the smallest subtorus having a coset containing $X$ and a coset containing $Y$.
    Then
    \[
    \dim (X:Y) \leq \dim(X) + \dim(Y) - \dim(H).
    \]
\end{lemma}

\section{Hodges Map and the Gravity Divisor}\label{sec:hodges-gravity-divisor}

The goal of this section is to prove our main result, Theorem \ref{thm:main_gravity}, which gives the geometric interpretation of Hodges formula for the reduced MHV gravity scattering amplitude $\widetilde{M_n}$.
We will do this by introducing an intrinsic vector bundle construction on the two-step flag variety $\F_n$.

\subsection{Hodges Matrix and Gravity Scattering Amplitudes}

We start by recalling the precise definitions of $\widetilde{M_n}$ and the Hodges matrix.
For an arbitrary spinor $z\in \C^2$ (not necessarily one of the $\lambda_i$'s from Section \ref{sec:preliminaries}) and an index $i$, we keep the notation $\langle j z\rangle$ to avoid confusion with the Pl\"ucker coordinates $p_{ij}$'s. 
For $i\in [n]$ let $\lambda_i \in \mathbb C^2 $ and $\widetilde\lambda_i \in \mathbb C^2$ be spinors satisfying \emph{momentum conservation}: for every fixed index $i$ and every spinor $z\in \mathbb C^2$,
\begin{equation}\label{eq:momentum-conservation-component}
    \sum_{j=1}^n q_{ij}\langle j z\rangle =0.
\end{equation}
Then we can define Hodges matrix.

\begin{definition}\label{def:Hodges_mat}
    Choose two auxiliary spinors $x,y\in \C^2$ such that $\langle i x\rangle \langle i y\rangle \neq 0$ for all $i$.
    The {\it Hodges matrix} is the symmetric $n\times n$ matrix $\Phi=\Phi(x,y)$ with  entries
    \[
      \Phi_{ij}\coloneqq \frac{q_{ij}}{p_{ij}},
      \mbox{ for } i\neq j,
    \]
    and
    \begin{equation}\label{eq:hodges-diagonal}
      \Phi_{ii}
      \coloneqq 
      -\sum_{j\neq i}
      \frac{q_{ij}}{p_{ij}}
      \frac{\langle jx\rangle \langle jy\rangle}
           {\langle i x\rangle\langle i y\rangle }.
    \end{equation}
\end{definition}

Finally, we can introduce the scattering amplitude of our main interest.

\begin{definition}\label{def:mhv_ampl}
    The \emph{reduced $n$-point MHV gravity amplitude} $\widetilde{M}_n$ in $\mathcal{N}=8$ supergravity (SUGRA) is the rational section on the spinor helicity variety $\mathcal{F}_n$ given by
    \[
        \widetilde{M}_n
        \coloneqq \frac{\det \Phi^{rst}_{ijk}}{p_{ij}p_{ik}p_{jk}p_{rs}p_{rt}p_{st}},
    \]
    where $\Phi^{rst}_{ijk}$ is the $(n{-3})$ square matrix obtained from the Hodges matrix $\Phi$  by deleting the rows $i,j,k$ and columns $r,s,t$.
\end{definition}

\begin{example}
    The \emph{$5$-point scattering amplitude} $\widetilde{M}_5$ is obtained from the $5\times 5$ Hodges matrix: 
    \[
    \Phi_5=
    \begin{bmatrix}
        \Phi_{11} & \frac{q_{12}}{p_{12}} & \frac{q_{13}}{p_{13}} & \frac{q_{14}}{p_{14}} & \frac{q_{15}}{p_{15}} \\
        \frac{q_{12}}{p_{12}} & \Phi_{22} & \frac{q_{23}}{p_{23}} & \frac{q_{24}}{p_{24}} & \frac{q_{25}}{p_{25}} \\
        \frac{q_{13}}{p_{13}} & \frac{q_{23}}{p_{23}} & \Phi_{33} & \frac{q_{34}}{p_{34}} & \frac{q_{35}}{p_{35}} \\
        \frac{q_{14}}{p_{14}} & \frac{q_{24}}{p_{24}} & \frac{q_{34}}{p_{34}} & \Phi_{44} & \frac{q_{45}}{p_{45}} \\
        \frac{q_{15}}{p_{15}} & \frac{q_{25}}{p_{25}} & \frac{q_{35}}{p_{35}} & \frac{q_{45}}{p_{45}} & \Phi_{55}
    \end{bmatrix}.
    \]
    We obtain
    \begin{align}\label{eq:A_5}
    \widetilde{M}_5 = 
        \frac{-p_{14} p_{25} q_{24} q_{15} + p_{24} p_{15} q_{14} q_{25}}{p_{12} p_{13} p_{23} p_{14} p_{24} p_{34} p_{15} p_{25} p_{35} p_{45}}.
    \end{align}
\end{example}

\subsection{Hodges bundle on the flag variety and quadrics}

In this subsection, we introduce the construction intrinsic to the partial flag variety $\F_n$, formulated in the language of vector bundles, which allows us to give a geometric interpretation of the reduced MHV scattering amplitude $\widetilde{M_n}$.

\medskip

Let $\mathcal{O}_{\F_n}$ be the structure sheaf on $\F_n$, and let $\cS\subset V\otimes\cO_{\F_n}$, with $V\cong \C^n$ be the rank-two tautological subbundle, whose fiber over $(L,W)\in \mathcal F_n$ is $L$.  
\begin{definition}
    For $i\in [n]$, define $s_i\in H^0(\mathcal F_n,\cS^\vee)$ as the restriction of the $i$-th coordinate functional on $\C^n$ to the 2-plane $L$. 
\end{definition}

\begin{remark}
    Choosing a basis of $L$, we may write it as span of the rows of a $2\times n$ matrix
    \[
        \Lambda=\begin{pmatrix}
        u_1&u_2&\cdots&u_n\\
        v_1&v_2&\cdots&v_n
        \end{pmatrix}.
    \]
    In these coordinates, $s_i=u_ie_1^\vee+v_ie_2^\vee\in L^\vee$, where $e_1$ and $e_2$ are elements from the chosen basis of $L$. 
    Thus $s_i$ can then be identified with the column spinor $\lambda_i=\begin{pmatrix}u_i\\v_i\end{pmatrix}$. 
\end{remark}

Note that the angle bracket  $\langle ij\rangle=p_{ij} = s_i\wedge s_j\in H^0(\mathcal F_n,\det\cS^\vee)$. 
This motivates the following definition.

\begin{definition}
    Define the line bundle $\cO_p(1)\coloneqq \det\cS^\vee$, whose sections are homogeneous polynomials of degree $1$ in the Pl\"ucker coordinates $p_{ij}$, analogously for $\cO_q(1)$.
    Then, define the line bundle $\cL:=\cO_q(1)\otimes \cO_p(-1)$, such that $z_{ij}=\frac{q_{ij}}{p_{ij}}$ are rational sections of $\cL$.
    On the open set $\F_n^\circ$ of $\mathcal F_n$ where $p_{ij}\neq 0$, they are regular sections of $\cL$.
\end{definition}

In the following, we shall restrict on an open locus of $\mathcal F_n$, the {\it Hodges locus}. This is the set $\mathcal F_n^\circ:=\Set{(L,W)\in \F_n \mid p_{ij}\neq 0\text{ for all }i<j}$.
On the Hodges locus, all off-diagonal entries $z_{ij}$ of the Hodges matrix are regular sections of the line bundle $\cL$.

\begin{definition}[Quadratic evaluation map]\label{def: quadric eval}
There is a natural vector bundle morphism
\[
 \ev_2:\Sym^2\cS\longrightarrow \cO_{\mathcal F_n}^{\oplus n}
\]
defined fiberwise by
\[
Q\longmapsto \bigl(Q(s_1),\ldots,Q(s_n)\bigr).
\]
Here $Q\in \Sym^2 L$ is a quadric on $L^\vee$; recall that $s_i\in L^\vee$. Denote its image by 
\[
    \cK \coloneqq \im(\ev_2)\subset \cO_{\mathcal F_n}^{\oplus n},
\]
where this image has constant rank.
\end{definition}

\begin{lemma}[Rank three on the Hodges locus]
    Let $n\geq 3$. On the Hodges locus $\mathcal F_n^\circ$, the evaluation map $\ev_2:\Sym^2\mathcal{S}\rightarrow \cO_{\F_n}^{\oplus n}$ is injective.
    Hence $\cK\simeq \Sym^2\cS$ is a rank-three subbundle of  $\cO_{\mathcal F_n^\circ}^{\oplus n}$.
\end{lemma}
\begin{proof}
    The assertion can be proven on each fiber.  Fix $(L,W)\in \mathcal F_n^\circ$ and let $Q\in \Sym^2 L$ be a quadric satisfying $Q(s_i)=0$ for all $1\leq i\leq n$. Since $p_{ij}=s_i\wedge s_j\neq0$ for all $i<j$, the points  $[s_1],\ldots,[s_n]\in \PP(L^\vee)\simeq \PP^1$ must be pairwise distinct. Therefore, since $n\geq 3$ and $Q$ is a binary quadric, the vanishing at all the $[s_i]$ implies that $Q$ vanishes identically.
    Hence $\ev_2$ is injective on every fiber of $\mathcal F_n^\circ$.  Clearly, $\Sym^2\cS$ has rank three, and so does its image.
\end{proof}

The usual Hodges diagonal is written using two reference spinors $x,y$. These two spinors give a quadric $f_{xy}$ in $\Sym^2 \cS$ defined by 
$f_{xy}(s) = x\wedge s \otimes y\wedge s\in (\det\cS^\vee)^{\otimes 2}$, where $s$ is a section of $\mathcal{S}^\vee$.
Intrinsically, this means choosing a $(\det\cS^\vee)^{\otimes 2}$-valued quadric on $\cS^\vee$.  

\begin{definition}[Hodges bundle]\label{def:Hodges_bundle}
    We call the rank $3$ bundle $\cE =\Sym^2\cS\otimes (\det\cS^\vee)^{\otimes 2}$ on $\mathcal F_n$ {\it Hodges bundle}. 
    A section $f$ of $\cE$ is then a quadric on $\cS^\vee$ with values in the line bundle $(\det\cS^\vee)^{\otimes 2}$. For such an $f$, define $U_f =\lbrace (L,W)\ | \ f(s_i)\neq 0\text{ for all } i=1,\ldots,n\rbrace\subset \mathcal F_n$.
\end{definition}

\begin{lemma}[The opens $U_f$ cover the Hodges locus]
    Hodges bundles $\mathcal{H}$ is globally generated.
    Moreover, if we denote by 
    \[
        \mathcal C_n^{\circ} \coloneqq \lbrace (L,W)\in \mathcal F_n\ | \ 0\neq s_i\in \det L^\vee \text{ for all }i\rbrace,
    \]
    then the open sets $U_f$, as $f$ varies over the global sections $H^0(\F_n,\cE)$ of $\cE$, cover $\mathcal C_n^{\circ}$.
    In particular, they cover $\mathcal F_n^{\circ}$.
\end{lemma}

\begin{proof}
    First observe that, because $\cS$ has rank two, there is a canonical isomorphism
    \[
     \cS\simeq \cS^\vee\otimes \det\cS.
    \]
    Indeed, fiberwise, for a 2-dimensional vector space $L$, the map is 
    \[
        L\longrightarrow L^\vee\otimes (L\wedge L), \ \ v\longmapsto \bigl(w\mapsto v\wedge w\bigr),
    \]
    which is an isomorphism.
    Taking symmetric squares on both sides gives $\Sym^2\cS\simeq\Sym^2\cS^\vee\otimes (\det\cS)^{\otimes2}$. 
    Tensoring by $(\det\cS^\vee)^{\otimes2}$, we obtain $\cE =\Sym^2\cS\otimes(\det\cS^\vee)^{\otimes2}\simeq \Sym^2\cS^\vee$. 
    Since $\cS\subset V\otimes\cO_{\mathcal F_n}$, dualizing gives a surjection $V^\vee\otimes\cO_{\mathcal F_n}\twoheadrightarrow \cS^\vee$ and so $\cS^\vee$ is globally generated. Thus $\Sym^2\cS^\vee\simeq\cE$ is globally generated too.
    This shows the first claim.
    
    Next, we show that the $U_f$ cover $\mathcal{C}^\circ_n$.
    Fix a point $(L,W)\in \mathcal C_n^\circ$.
    By definition of $\mathcal C_n^\circ$, the sections $ s_1,\ldots,s_n\in L^\vee$ are all nonzero, and therefore determine finitely many points $[s_1],\ldots,[s_n]\in \mathbb P(L^\vee)\simeq \mathbb P^1$.
    Recall that $\operatorname{Sym}^2 L = \operatorname{Sym}^2\big((L^\vee)^\vee\big)$ is naturally the vector space of homogeneous quadratic forms on
    $L^\vee$.
    For every $i$, the condition that a quadratic form $Q\in\operatorname{Sym}^2L$ vanishes at $s_i$, i.e. $Q(s_i)=0$,
    defines a two-plane in $\operatorname{Sym}^2L$.
    Hence, there is a Zariski dense open set of elements $Q\in \Sym^2L$, such that $ Q(s_i)\neq0$ for every $i=1,\ldots,n$.
    Choose also a nonzero element $\omega\in(\det L^\vee)^{\otimes2}$.
    Then
    \[
        f_{(L,W)}:=Q\otimes\omega \in \operatorname{Sym}^2L\otimes(\det L^\vee)^{\otimes2} =\mathcal H_{(L,W)}
    \]
    satisfies  $f_{(L,W)}(s_i) = Q(s_i)\,\omega \neq0$, for every $i$. Here $f_{(L,W)}(s_i)\in(\det L^\vee)^{\otimes2}$
    is the value of the quadratic form $f_{(L,W)}$ at the vector
    $s_i\in L^\vee$.
    Since $\mathcal H$ is globally generated, the evaluation map $H^0(\mathcal F_n,\mathcal H) \longrightarrow \mathcal H_{(L,W)}$ is surjective.
    Hence there exists a global section $f\in H^0(\mathcal F_n,\mathcal H)$ such that $f(L,W)=f_{(L,W)}$.
    Consequently, we have $f(L,W)(s_i)\neq0 \qquad\text{for every }i$, 
    and therefore $(L,W)\in U_f$.
    So the $U_f$'s cover $\mathcal C_n^{\circ}$.
    Finally, it is clear that $\mathcal C_n^{\circ}\supset \mathcal F_n^{\circ}$, because $p_{ij}\neq 0$ implies that $s_i, s_j\neq 0$. 
\end{proof}

\begin{remark}
    Note that using the sections $s_k\in H^0(\F_n,\mathcal{S}^{\vee})$, the momentum-conservation equations can be written  as
    \[
     \sum_{k=1}^n q_{ik}s_k=0
    \qquad\text{for every }i=1,\ldots,n.
    \]
    In coordinates this is the identity $\sum_{k=1}^n q_{ik} \lambda_k=0$. Indeed, as we mentioned in \S\ref{sec:preliminaries}, if $\lambda_k$ denotes the $k$-th angle spinor column and $\tilde\lambda_k$ the square spinor, then momentum conservation is $\sum_{k=1}^n \lambda_k\tilde\lambda_k=0$.
    Contracting with $\tilde\lambda_i$ on the square-spinor side gives $\sum_{k=1}^n q_{ik} \lambda_k=0$.
\end{remark}

\medskip

Let $f$ be a section of $\cE$. We now define a generalized Hodges matrix. 

\begin{definition} \label{def:generalized_Hodges_mat}
    On $\mathcal F_n^\circ\cap U_f$, define the $n\times n$ symmetric matrix $\Phi^{(f)}$, called \emph{generalized Hodges matrix}, with entries in the line bundle $\cL$ as follows:
    \[
     \Phi^{(f)}_{ij} 
     \coloneqq z_{ij}=\frac{q_{ij}}{p_{ij}}, \mbox{ and } 
     \Phi^{(f)}_{ii} \coloneqq -\frac{\sum_{k\neq i} z_{ik}f(s_k)}{f(s_i)}.
    \]
\end{definition}

\begin{remark}
    Firstly, we will show that this definition is independent of the choice of section $f$ in Proposition \ref{prop:indepedence_quadric_f}.
    Also note that this quotient makes sense because $f(s_k)$ and $f(s_i)$ are sections of the same line bundle $(\det \ms^{\vee})^{\otimes 2}$, so $f(s_k)/f(s_i)$ is a scalar regular function on $U_f$. Hence $\Phi^{(f)}_{ii}$ is still a section of $\cL$.
    Moreover, if $f=f_{x,y}$ is the quadric $f_{x,y}(s)= s\wedge x\otimes  s\wedge y$, then the diagonal entries of $\Phi(f_{x,y})$ are those of the usual Hodges matrix. 
\end{remark}

\begin{remark}
Given the global sections $s_1,\ldots, s_n\in H^0(\mathcal F_n, \mathcal S^{\vee})$, we have that every global section of the Hodges bundle $\mathcal H$ is of the form 
\[
f = \sum_{1\leq i\leq j\leq n} c_{ij}f_{ij}, \ \ c_{ij}\in \C,
\]
where $f_{ij}$ are quadrics such that 
$f_{ij}(s)=(s_i\wedge s)(s_j\wedge s)$. Note that these are $\det \mathcal S$-valued. Thus $(c_{ij})$ is a complex $n\times n$ symmetric matrix. When $f=f_{x,y}$ as in the usual Hodges matrix, the matrix $(c_{ij})$ has rank one. For general global sections this matrix of coefficient is of course full-rank. 
\end{remark}

We will need the following technical lemma in order to prove the independence of the choice of quadric $f\in H^0(\F_n,\mathcal{H})$ in this definition.

\begin{lemma}\label{lem: vanishing}
    Fix $i$ and let $h$ be any local section of $\Sym^2\cS\otimes \mathcal E$, for some line bundle $\mathcal E$ on $\F_n$. 
    Suppose that, at a given point of $\mathcal F_n^\circ$, $h(s_i)=0$. Then, at that point, we also have $\sum_{k\neq i}z_{ik}h(s_k)=0$.
\end{lemma}
\begin{proof}
    The statement is local, so we may trivialize all the line bundles involved and choose a basis of the rank-two vector space $L$.  On such a local chart, $h$ is a quadric on the two-dimensional vector space $L^\vee$.  The condition $h(s_i)=0$ says that the quadratic polynomial $h$ vanishes at the point $[s_i]\in \PP(L^\vee)$. A quadric vanishing at $[s_i]$ is divisible by the linear form $s\longmapsto s\wedge s_i$. Hence there exists a linear form $\ell$ such that $h(s)=(s\wedge s_i)\ell(s)$. Evaluating at $s_k$, we obtain $h(s_k)=(s_k\wedge s_i)\ell(s_k)$. Now, $s_k\wedge s_i=-s_i\wedge s_k=-p_{ik}$. So $h(s_k)=-p_{ik}\ell(s_k)$.
    Therefore
    \[
    \sum_{k\neq i}z_{ik}h(s_k) = \sum_{k\neq i}\frac{q_{ik}}{p_{ik}}\bigl(-p_{ik}\ell(s_k)\bigr)= -\sum_{k\neq i}q_{ik}\ell(s_k) = -\sum_{k=1}^n q_{ik}\ell(s_k),
    \]
    where, in the last equality, we added the term for $k=i$ which does not contribute. By momentum conservation, $\sum_{k=1}^n q_{ik}s_k=0$.
    Pairing this identity with the linear form $\ell$, we obtain that $\sum_{k=1}^n q_{ik}\ell(s_k)=0$. 
    Thus $\sum_{k\neq i}z_{ik}h(s_k)=0$. 
\end{proof}

\begin{proposition}[Independence of the quadric]\label{prop:indepedence_quadric_f}
    Let $f$ and $g$ be two local sections of $\cE$.  On $\mathcal F_n^\circ\cap U_f\cap U_g$, one has $\Phi^{(f)}=\Phi^{(g)}$. 
    Consequently, the local matrices $\Phi^{(f)}$ glue to a canonical matrix
    \[
    \Phi\in H^0\bigl(\mathcal F_n^\circ,\Sym^2(\cO_{\mathcal F_n^\circ}^{\oplus n})^\vee\otimes\cL\bigr).
    \]
    Equivalently, $\Phi$ gives a rational projective map, the \emph{Hodges map}: 
    \[
    H: \mathcal F_n \dashrightarrow \P(\Sym_\C(n)), \ \ (L,W)\longmapsto [\Phi(L,W)].
    \]
\end{proposition}

\begin{proof}
    The off-diagonal entries are the same by definition. We only need to compare the diagonal entries. Fix $i$.  After trivializing $\cE$ locally, regard $f(s_i)$ and $g(s_i)$ as nonzero scalars.  Define $h=g(s_i)f-f(s_i)g$. Then $h(s_i) = 0$. By Lemma \ref{lem: vanishing}, $\sum_{k\neq i}z_{ik}h(s_k)=0$.
    Expanding this identity gives
    \[
     g(s_i)\sum_{k\neq i}z_{ik}f(s_k)-f(s_i)\sum_{k\neq i}z_{ik}g(s_k)=0.
    \]
    Dividing by $f(s_i)g(s_i)$ yields
    \[
     \frac{\sum_{k\neq i}z_{ik}f(s_k)}{f(s_i)}=
    \frac{\sum_{k\neq i}z_{ik}g(s_k)}{g(s_i)}.
    \]
    Hence $\Phi^{(f)}_{ii}=\Phi^{(g)}_{ii}$. This proves the independence of $f$.  Since the open sets $U_f$ cover $\mathcal F_n^\circ$, the local matrices glue. Finally, the entries of $\Phi$ are sections of the same line bundle $\cL$.  Hence, after trivializing $\cL$, they give an ordinary symmetric matrix; changing the trivialization scales all entries by the same nonzero function.  Therefore the projective point $[\Phi]\in\PP(\Sym_n)$ is well-defined.  This describes the Hodges rational map $H$.
\end{proof}

\subsection{The Kernel of the Generalized Hodges Matrix}

It was already shown in Hodges' original work \cite{Hodges2012}, that the Hodges matrix from Definition \ref{def:Hodges_mat} has rank $n-3$.
Motivated by this observation here we will study the kernel of the generalized Hodges matrix from Definition \ref{def:generalized_Hodges_mat} in the global language of line bundles.
This will allow us to derive the global version of Hodges formula in Definition \ref{def:mhv_ampl}.
The main result can be stated as follows.

\begin{proposition}[The Hodges kernel]\label{prop:Hodge_kernel}
    On the Hodges locus $\mathcal F_n^\circ$, the rank-three bundle $\cK$
    is contained in the kernel of $\Phi$. Equivalently, for every local quadric $Q\in\Sym^2\cS$, the vector $v_Q =\bigl(Q(s_1),\ldots,Q(s_n)\bigr)$ satisfies
    $\Phi v_Q=0$. 
\end{proposition}

\begin{proof}
    By Proposition \ref{prop:indepedence_quadric_f} it suffices to work on an open set $U_f$, where we can use the local formula $\Phi^{(f)}$ for $\Phi$.
    Fix $i\in [n]$.
    We have to show that $\Phi_{ii}Q(s_i)+\sum_{k\neq i}z_{ik}Q(s_k)=0$.
    By definition of the generalized Hodges matrix as in Definition \ref{def:generalized_Hodges_mat}, we locally have 
    \[
    \begin{aligned}
            \Phi_{ii}Q(s_i)+\sum_{k\neq i}z_{ik}Q(s_k)
            &=-\frac{Q(s_i)}{f(s_i)}\sum_{k\neq i}z_{ik}f(s_k)
              +\sum_{k\neq i}z_{ik}Q(s_k)\\
            &=\frac{1}{f(s_i)}\sum_{k\neq i}z_{ik}
            \bigl(-Q(s_i)f(s_k) + f(s_i)Q(s_k)\bigr).
    \end{aligned}
    \]
    Consider the quadric $h=f(s_i)Q-Q(s_i)f\in \mathrm{Sym}^2\mathcal{S}\otimes (\det \mathcal S^{\vee})^{\otimes 2}$. Then $h(s_i)=0$. Hence $\sum_{k\neq i}z_{ik}h(s_k)=0$ by Lemma \ref{lem: vanishing}. Therefore $\Phi_{ii}Q(s_i)+\sum_{k\neq i}z_{ik}Q(s_k)=0$. Since this holds for every row $i$, we derive $\Phi v_Q=0$.
\end{proof}

\begin{corollary}\label{cor: blowup}
    The Hodges rational map induces a rational map
    \[
    \widetilde H:\F_n\dashrightarrow \Set{([A],K)\in\PP(\Sym_n)\times\Gr(3,n):K\subseteq\ker A},
    \]
    on the Hodges locus $\mathcal F_n^\circ$ defined by
    \[
    (L,W)\longmapsto \bigl([\Phi(L,W)],\cK_{(L,W)}\bigr),
    \]
    where $\mathcal{K}_{(L,W)}$ is the fiber of $\mathcal{K}$ at $(L,W)$.
    Moreover, the incidence space on the right is isomorphic to the blowup of the variety of symmetric matrices of rank at most $(n-3)$ along the subvariety of those whose rank is at most $(n-4)$. 
\end{corollary}

\subsection{The Induced Determinant and the MHV gravity divisor}

In this paragraph we finally give the geometric interpretation of the MHV gravity amplitude.
We start, however, with a technical result.

\begin{lemma}\label{lem: induced det}
    Let $A$ be an $n\times n$ symmetric complex matrix and suppose that we are given a $3$-dimensional subspace $K\subseteq \ker(A)$. Let $Q=\C^n/K$.
    \begin{itemize}
        \item[(i)] Then $A$ descends to a symmetric linear map
                \[
                \overline A: Q \longrightarrow Q^\vee.
                \]
        \item[(ii)] Choosing the standard volume element on $\C^n$, the \emph{induced determinant} $\det \overline{A}\in (\det Q^\vee)^{\otimes 2}$ can be identified with an element $\overline{\Delta}\in (\det K)^{\otimes 2}$.
        \item[(iii)] Denote by $A_R^C$ the $(n-3)\times(n-3)$ matrix obtained from $A$ by deleting the rows indexed by $R$ and the columns indexed by $C$.
        Then,
        \[
            \det A_R^C = \pm \langle \kappa_R\otimes \kappa_C\,,\, \overline \Delta\,\rangle \in \C,
        \]
        where $\kappa_R$ and $\kappa_C$ in $(\det K)^\vee$ are the Pl\"ucker coordinate functionals of $K$ labeled by $R$ and $C$, respectively, and $\langle-,-\rangle$ is the natural pairing $(\det K^\vee)^{\otimes2}\times (\det K)^{\otimes 2}\to \C$.
        In particular, for any $R$ and $C$ the minor $\det A^C_R$ is obtained by evaluating the \emph{induced determinant} $\det(\overline A)$ on the functionals $\kappa_R$ and $\kappa_C$.
    \end{itemize}
\end{lemma}
\begin{proof}
    We start by showing (i). 
    Because $K\subseteq \ker(A)$, the bilinear form associated to $A$ vanishes whenever one argument lies in $K$.  Let $Q=\C^n/K$. Hence it induces a well-defined bilinear form on the quotient $Q$, equivalently a symmetric map
$\overline A: Q\to Q^\vee$. 
In the rest of the proof, we implicitly apply the canonical identification between $(\det Q)^{\vee}$ and $\det Q^{\vee}$. 

We show statement (ii). Choose the standard volume form $\omega =e_1\wedge\ldots\wedge e_n$ on $\C^n$, then $\det Q^\vee \cong \det K$, as can be seen by the short exact sequence
\begin{align}\label{eq:SES_maps}
        0\to K\to \C^n\to Q\to 0.
    \end{align}
Now recall that the determinant of a map $Q\to Q^\vee$ is an element of
$\det Q^\vee \otimes \det Q^\vee \cong (\det K) ^{\otimes 2}$. Then, $\overline{\Delta}\in (\det K)^{\otimes 2}$ is the image of $\det \overline{A}$ under this isomorphism.

It remains to show (iii).  
Let $U=(u_1,u_2,u_3)$ be the $n\times 3$ matrix spanning $K$, and let $Z=(z_1,\ldots,z_{n-3})$ be a $n\times (n-3)$ matrix such that the $n\times n$ matrix $B=[UZ]$ has non-vanishing determinant.
Let $z=\bar z_1\wedge\cdots\wedge\bar z_{n-3}\in \det Q$, and take its dual $z^\vee\in \det Q^\vee$ with $\langle z^\vee,z\rangle=1$. Then, $\overline{A}$ has a matrix representative $Z^TAZ$ (with respect to the basis given by the $z_i$).
    We obtain
    \begin{align}\label{eq:detA_1}
    \det(\overline A)=\det(Z^TAZ)(z^\vee)^{\otimes2}\in (\det Q^\vee)^{\otimes 2}.
    \end{align}
    The determinant isomorphism from the short exact sequence \eqref{eq:SES_maps}
    \[
    \mu:\det K\otimes\det Q\cong \det \C^n
    \]
    satisfies $\mu(u_1\wedge u_2 \wedge u_3\otimes z)=\det(B)\omega$.
    Thus, the identification induced by $\omega$ gives a canonical isomorphism  
    $G_{\omega}: \det K\rightarrow \det Q^{\vee}$ defined by $G_{\omega}(k)(q)\omega = \mu(k\otimes q)$. One can check that its inverse is the isormorphism
\[
F_\omega:\det Q^\vee\longrightarrow\det K,
\qquad
\lambda\longmapsto
(\operatorname{id}_{\det K}\otimes\lambda)(\mu^{-1}(\omega)).
\]
This linear map is such that  
\begin{align}\label{eq:detA_2}
 z^\vee\mapsto (\operatorname{id}_{\det K}\otimes z^\vee)(\mu^{-1}(\omega))
 &=\frac{1}{\det B}\langle \operatorname{id}_{\det K}, u_1 \wedge u_2 \wedge u_3\rangle \otimes \langle z^\vee ,z \rangle \\
&=\frac{u_1 \wedge u_2 \wedge u_3}{\det B},
\end{align}
since $\mu^{-1}(\omega)=(u_1 \wedge u_2 \wedge u_3\otimes z)/\det B$ and $\langle z^\vee, z\rangle=1$.
Combining \eqref{eq:detA_1} and \eqref{eq:detA_2} yields
\begin{align}\label{eq:detDelta}
\overline\Delta
\coloneqq F_\omega^{\otimes2}\bigl(\det(\overline A)\bigr)=\frac{\det Z^TAZ}{(\det B)^2}\,(u_1 \wedge u_2 \wedge u_3)^{\otimes2}.
\end{align}
Moreover, as $AU=0$ we have
\[
B^TAB = \begin{pmatrix}
0_{3\times 3}& 0_{3\times (n-3)} \\
0_{(n-3)\times 3} & Z^TAZ \\
\end{pmatrix}.
\]
Therefore, 
$A = B^{-T}\begin{pmatrix}
0_{3\times 3}& 0_{3\times (n-3)} \\
0_{(n-3)\times 3} & Z^TAZ \\
\end{pmatrix} B^{-1}$.
Recall that for a matrix $A$ we denote by $A^C_R$ the matrix obtained from $A$ by {\it deleting} the rows and columns indexed by $R$ and $C$, respectively.
Then, by applying Cauchy--Binet formula twice to the matrix equation above, we get
    \begin{align}
        \det A_R^C
        &= \sum_{\substack{I,J\subseteq\{1,\ldots,n\}\\ |I|=|J|=3}} \det (B^{-T})_R^I\, \det (B^TAB)_I^J\, \det (B^{-1})_J^C\\
        &= \det (B^{-T})_R^{S^c}\, \det(Z^TAZ)\, \det (B^{-1})_{S^c}^C.
    \end{align}
    where $S=\{4,\ldots,n\}$, and $R^c,C^c$ denote the complements of $R,C$ in $[n]$. The last equality follows because the first three rows and columns of $B^TAB$ vanish, so only $I=J=S^c$ can contribute.
    Then, using Jacobi's identity on complementary minors,
    \begin{align}
        \det (B^{-T})^{S^c}_R &= \pm \frac{\det (B^T)^{R^c}_S}{\det B} =\pm\frac{\det B^{R^c}_S}{\det B},\\
        \det (B^{-1})_{S^c}^{C} &=\pm \frac{\det B_{C^c}^{S}}{\det B}.
    \end{align}
Note that $B_{R^c}^{S}=[U Z]_{R^c}^{S}=U_{R^c}$ and $B_{C^c}^{S}=[UZ]_{C^c}^{S}=U_{C^c}$; it follows that
    \[
        \det A_R^C = \pm \frac{\det U_{R^c}\,\det U_{C^c}}{\det(B)^2}\det(Z^TAZ).
    \]
    To establish (iii), we now compute the pairing $\langle\kappa_R\otimes\kappa_C,\overline\Delta\rangle$.
    For $R=\{r_1<r_2<r_3\}$ let 
    \[
        \kappa_R= \left.(e_{r_1}^\vee\wedge e_{r_2}^\vee\wedge e_{r_3}^\vee)\right|_{\det K} \in \det K^\vee.
    \]
    and analogously for $\kappa_C$.
    Then, $\kappa_R(u_1 \wedge u_2 \wedge u_3)=\det U_{R^c}$ and $\kappa_C(u_1 \wedge u_2 \wedge u_3)=\det U_{C^c}$.
    Therefore,
    \begin{align}
        \langle\kappa_R\otimes\kappa_C,\overline\Delta\rangle
        &=\frac{\det(Z^TAZ)}{(\det B)^2}\,
          \kappa_R(u_1 \wedge u_2 \wedge u_3)\kappa_C(u_1 \wedge u_2 \wedge u_3)\\
        &=\frac{\det(Z^TAZ)}{(\det B)^2}\,
          \det U_{R^c}\,\det U_{C^c}=\pm\det A^C_{R}. \qedhere
    \end{align}
\end{proof}

Now, on the Hodges locus $\mathcal F_n^\circ$, let $\cQ:=\cO_{\mathcal F_n^\circ}^{\oplus n}/\cK$. This is a vector bundle of rank $n-3$. 
Since $\Phi$ is symmetric and $\mathcal{K}$ is in the kernel of $\Phi$ by Proposition \ref{prop:Hodge_kernel}, it descends to a symmetric $\cL$-valued form on $\cQ$, i.e. $\overline\Phi:\cQ\longrightarrow \cQ^\vee\otimes\cL$ by Lemma \ref{lem: induced det}.
Equivalently, it is a section of $\Sym^2\cQ^\vee\otimes\cL$. 
Recall that the determinant of a map $\cQ\to\cQ^\vee\otimes\cL$ is a (global) section of 
\[
    \det(\cQ^\vee\otimes\cL)\otimes(\det\cQ)^{\vee}
    =\det\cQ^\vee\otimes\cL^{\otimes (n-3)}\otimes\det\cQ^\vee
    =(\det\cQ^\vee)^{\otimes2}\otimes\cL^{\otimes (n-3)}.
\]
Thus, $\Delta_H \coloneqq \det(\overline\Phi)\in H^0\left(\F_n^\circ, (\det\cQ^\vee)^{\otimes2}\otimes\cL^{\otimes(n-3)}\right)$.
This motivates the following definition.

\begin{definition}\label{def:gravity_divisor}
     
    Denote the vanishing locus of the regular section $\Delta_H$ on Hodges locus $\F_n^\circ$ by $D_H^\circ\coloneqq \mathbb{V}(\Delta_H)\subset \F_n^\circ$, which we call the \emph{open gravity divisor}.
    Its closure $D_H\subset \F_n$ in the flag variety $\F_n$ is called the \emph{gravity divisor}.
\end{definition}

\begin{remark}
It is known by \cite{Hodges2012} that $\Phi$ has generic rank $n-3$ on $\mathcal F_n^\circ$.  
\end{remark}

We can now record the following version of our main Theorem \ref{thm:main_gravity}, using the language of vector bundles introduced above.

\begin{theorem}[The MHV gravity divisor] \label{thm:gravity_divisor}\hfill 
    Let $n \geq 5$ and $\F_n^\circ$ be the Hodges locus in $\F_n$ where all $p_{ij}\neq 0$.
    \begin{enumerate}[itemsep = 0pt]
        \item[(i)] The open gravity divisor $D_H^\circ=\mathbb{V}(\Delta_H)\subset \mathcal F_n^\circ$ is an effective Cartier divisor. 
        Set-theoretically it is the locus of points $(L,W)\in \F_n^\circ$ where the induced bi-linear form $\overline\Phi:\C^n/\mathcal{K}_{(L,W)}\to (\C^n/\mathcal{K}_{(L,W)})^\vee\otimes\cL_{(L,W)}$ fails to be invertible, with $\mathcal{K}_{(L,W)}$ being the fiber of $\mathcal{K}$ over $(L,W)$.
        In symbols we have
        \[
            D_H^\circ = \left\{ (L,W)\in \F_n^\circ\mid \det(\overline{\Phi}(L,W))=0\right\}\subset \F_n^\circ.
        \]
        
        \item[(ii)] For $R = \lbrace i<j<k\rbrace$ and $C=\lbrace r<s<t\rbrace$, we have
        \[
            \Delta_H= \pm \frac{\det\Phi^C_R}{p_{ij}p_{ik}p_{jk}\,p_{rs}p_{rt}p_{st}}\in H^0\left(\F_n^\circ, (\det\cQ^\vee)^{\otimes2}\otimes\cL^{\otimes(n-3)}\right),
        \]
        so it coincides with the MHV amplitude $\widetilde{M}_n$ for $n$ particles in $\mathcal N=8$ supergravity.
        \item[(iii)]  The pullback under the Hodges rational map $\wt{H}$ of the exceptional divisor $E$ of the blowup of the variety of symmetric $n$ matrices of rank at most $n-3$ in its singular locus gives the open gravity divisor $D_H^\circ$.
        Put in symbols we have
        \[
            \wt{H}^{\,*}E\cong D_H^\circ.
        \]
    \end{enumerate}
\end{theorem}
\begin{proof}
    We start with (i). This follows from the discussion after Lemma \ref{lem: induced det}.
    Firstly, $\Delta_H$ is by definition a regular global section of the line bundle $(\det\cQ^\vee)^{\otimes2}\otimes\cL^{\otimes(n-3)}$ on the smooth variety $\F_n^\circ$ hence its vanishing locus is Cartier and effective.
    Since $\Delta_H=\det \overline{\Phi}$ it is precisely the locus of points in $(L,W)\in \F_n^\circ$ where $\overline{\Phi}(L,W)$ is not invertible.

Next we prove statement (ii). Choose a basis of $L$, then write $s_a=(u_a,v_a)\in \det L^\vee$, and choose a basis of $\Sym^2L$ as $(x^2,xy,y^2)$.
    
Then, by Proposition \ref{prop:Hodge_kernel} $\mathcal{K}$ is generated by evaluating $\Sym^2L$ on the $s_a$.
Hence, the $a$-th row of the $n\times 3$ matrix $U$ spanning (a subspace of) the kernel $\mathcal{K}_{(L,W)}$ is given by $(u_a^2,u_av_a,v_a^2)$.
Therefore, 
\[
\kappa_{C}=\det U_{C^{c}}=(u_rv_s-u_sv_r)(u_rv_t-u_tv_r)(u_sv_t-u_tv_s)=p_{rs}p_{rt}p_{st}.
\]
Similarly for $\kappa_R$. Then, the claim follows from Lemma \ref{lem: induced det}.
Moreover, by Proposition \ref{prop:indepedence_quadric_f} $\Phi$ is independent of the choice of quadrics, used to define the matrix locally, hence we can assume that $\Phi$ is locally the Hodges matrix. It remains to prove (iii): this follows directly from Corollary \ref{cor: blowup} together with (i). 
\end{proof}

\section{Pentad and Affine Slope Varieties}\label{sec:pentads}

Motivated by the construction of Hodges matrix and the Hodges rational map in the previous section, we study the Hadamard ratio of two Grassmannians $\Gr(2,n)$, which we will call the \emph{pentad variety $\pen_n$}.
The goal of this section is to study the dimension, defining equations, its tropicalization.
The results are summarized in Theorem \ref{thm:main_pentad}.
As mentioned in the introduction we will also consider the notion of little group scaling in our setup.
This will lead us to the affine slope variety Definition \ref{def:affine_slope} and Corollary \ref{cor:main_aff_to_pen} makes this connection precise.
The affine slope variety will also be a key ingredient in the proof of our main result about the pentad variety itself.

\subsection{Pentad Variety}

We start by introducing the pentad variety more formally.
Let $n\geq 3$.
The Grassmannian $\gr(2,n)\hookrightarrow\P(\wedge^2\C^n)$ is the subvariety of skew-symmetric $n\times n$ matrices $\P(\wedge^2\C^n)=\P^N$ with $N={\binom{n}{2}-1}$ given by matrices $Q=(q_{ij})$ such that $\rk Q=2$.
Take a second copy of $\gr(2,n)$ and denote the entries of its rank $2$ skew symmetric matrices by $p_{ij}$.
In this section we denote the polynomial ring generated by the $z_{ij}$ by $Z_n\coloneqq \C[z_{ij}\mid 1\leq i<j\leq n]$.
%and $Z_n$ the $n\times n$ matrix with entries $z_{ij}$; the diagonal entries are to be determined later.

\begin{definition}\label{def:pentad_variety}
    Denote the Hadamard ratio map for the two Grassmannians by
    \begin{align}\label{eq:Hadamard_ratio}
        \Psi_n:\gr(2,n)\times \gr(2,n) \ratto \PP^{N},\quad (q_{ij},p_{ij})\mapsto z_{ij}=\frac{q_{ij}}{p_{ij}}.
    \end{align}
    Then, the \emph{pentad variety} $\pen_n\subset \P^{N}$ is the Hadamard ratio $\pen_n\coloneqq \Gr(2,n):\Gr(2,n)$.
\end{definition}

\begin{remark}
    The rational Hadamard ratio map can be resolved to a morphism using the {\it De Concini-Procesi sequential blowup} construction \cite{DeConProc}. 
\end{remark}

\begin{definition}
    Let $\pen_n\subset \P^{N}$ be the closure of the image of the $n$-th Hadamard ratio map $\Psi_n$. This is called \emph{pentad variety}. 
\end{definition}

\begin{remark}\label{rem:pentad_exm}
    The Hadamard ratio map is dominant for $n \leq 4$. The first
    interesting case is $n=5$. Here we find that ${\rm Pen}_5$ is a quintic
    hypersurface in $\PP^9$. Namely, this hypersurface is defined by the polynomial
    \begin{align}
     \label{eq:pentad}
    \begin{matrix}
         \text{pen}_{12345}=&\quad z_{12} z_{13} z_{24} z_{35} z_{45}
        -z_{12} z_{13} z_{25} z_{34} z_{45}
        -z_{12} z_{14} z_{23} z_{35} z_{45}
        +z_{12} z_{14} z_{25} z_{34} z_{35} \\
        &+\,z_{12} z_{15} z_{23} z_{34} z_{45}
        -z_{12} z_{15} z_{24} z_{34} z_{35}
        +z_{13} z_{14} z_{23} z_{25} z_{45}
        -z_{13} z_{14} z_{24} z_{25} z_{35} \\
        &\,-\,z_{13} z_{15} z_{23} z_{24} z_{45}
        +z_{13} z_{15} z_{24} z_{25} z_{34}
        +z_{14} z_{15} z_{23} z_{24} z_{35}
        -z_{14} z_{15} z_{23} z_{25} z_{34}.
    \end{matrix}
    \quad
    \end{align}
    This is known as the {\em pentad} in algebraic statistics \cite{pentad}, which gives the name to the variety we are interested in. 
    Its $12$ monomials correspond to the $12$ oriented five-cycles of a pentagon.
\end{remark}

Suppose now that $n \geq 6$.
For any $5$-tuple $I$ in $[n]$, there is an analogous pentad, but now in the unknowns $z_{ij}$ where $i,j \in I$.
This is obtained by relabeling \eqref{eq:pentad}.
We denote this pentad by $\rm pen_{I}$.
Let $J^p_n$ denote the ideal in  the polynomial ring $Z_n$ which is generated by all $\binom{n}{5} $ pentads, i.e.
\begin{align}\label{eq:def_Jpn}
    J^p_n \coloneq \left( \operatorname{pen}_I \mid I\subset [n] \text{ s.t. }|I|=5\right) \subset Z_n.
\end{align}
By Remark \ref{rem:pentad_exm} the pentads vanish on the image of the Hadamard ratio map \eqref{eq:Hadamard_ratio}, so we have $\Pen_n \subseteq \mathbb{V}(J^p_n)$.
However, the ideal $J^p_n$ is not prime.
Since the base locus of the Hadamard ratio map $\Psi_n$ includes the Pl\"ucker hyperplane section of $\Gr(2,n)$, we consider the saturation
\begin{align}\label{eq:def_Ipn}
    I^p_n \,\,\coloneqq \,\, \biggl( \,J^p_n \,: \,\bigg[\prod_{1 \leq i < j \leq n} \!\!\! z_{ij} \,\bigg]^\infty \,\biggr)\subset Z_n.
\end{align}
This leads us to the following conjecture.
\begin{conjecture} \label{con:primeness}
    The ideal $I^p_n$ is prime of codimension $\binom{n-3}{2}$.
    It defines the irreducible variety
    $$ \Pen_n \, = \,\mathbb{V}(I^p_n). $$
    In words, the Hadamard ratio of the two Grassmannians is cut out by the pentad constraints up to saturation.
\end{conjecture}

\begin{remark}
    This conjecture can be proven for the case $n=6,7$ with a direct computation in \texttt{Macaulay2}.
    An interesting combinatorial question is to understand the combinatorics of the generators of $I^p_n$  that are not in $J^p_n$.
    We collect the generators of $I^p_n$ for $n=5,6,7,8$ in Table \ref{tab:betti_tables_I_n} below. (For $n=8$ we could only verify that there are no more generators of degree $8$.)
\begin{table}[H]
    \centering
    \caption{Number of generators of $I^p_n$ for small $n$ by degree.}
    \label{tab:betti_tables_I_n}
    \renewcommand{\arraystretch}{1.2}
    \begin{tabular}{llll|l}
        \toprule
         $\deg$ & $n=5$ & $n=6$ & $n=7$ & $n=8$ \\
        \midrule
         $5$
            & $1$ 
            & $6$ 
            & $21$
            & $56$\\
        $6$ 
            & $0$ 
            & $5=5\binom{6}{6}$ 
            & $35=5\binom{7}{6}$
            & $140 = 5\binom{8}{6}$\\
        $7$ 
            & $0$ 
            & $0$
            & $30 = 5\cdot 6\binom{7}{7}$ 
            & $240 = 5\cdot 6 \binom{8}{7}$ \\
        \bottomrule
    \end{tabular}
\end{table}
\end{remark}

For convenience, we recall our main result for the pentad variety, Theorem \ref{thm:main_pentad}, here again as follows.

\begin{theorem}\label{thm:pentad_variety}
    The pentad variety $\pen_n\subset \P^N$ with $N=\binom{n}{2}-1$ is irreducible of dimension $3n-7$, with prime vanishing ideal $\sqrt{I^p_n}$.
    So in particular,
    \[
        \mathbb{V}\left(\sqrt{I^p_n}\right) = \mathbb{V}(I^p_n) = \pen_n.
    \]
    Moreover, for $n\leq 10$ the ideal $I^p_n$ is radical, hence $I^p_n$ is the prime vanishing ideal. Equivalently, for $n\leq 10$ Conjecture \ref{con:primeness} holds.
\end{theorem}

In the remainder of this section we prove the statements in this theorem.
We start with the claim about the dimension, which is asserted by the following proposition.
Denote by $\ogr(2,n)$ the open Grassmannian defined by removing all Pl\"ucker hyperplane sections from $\gr(2,n)$ in its Pl\"ucker embedding.
Moreover, let $\mathring{\pen}_n$ be the image of the restriction of $\Psi_n$ to these open sets.

\begin{proposition}\label{prop:dimPen}
    The pentad variety $\pen_n\subset \P^N$ has dimension $3n-7$, equivalently it has codimension $\binom{n-3}{2}$ in $\P^N$ with $N={\binom{n}{2}-1}$.
\end{proposition}
\begin{proof}
    We start by giving an upper bound for the dimension.
    First note that $H=(\C^\times)^{n-1}\cong (\C^\times)^n/\C^\times$ is a subtorus of $(\C^\times)^{N}$ acting on $\ogr(2,n)$.
    Hence, 
    \begin{align}\label{eq:upper_bound_dim}
        \dim \mathring \pen_n &= \dim (\ogr(2,n):\ogr(2,n)) \\
        &\leq  2 \dim \ogr(2,n) - \dim (\C^\times)^{n-1} = 4n-8 - (n-1) = 3n-7
    \end{align}
    by Lemma \ref{lem:tropDim}. Next, we show that this bound is tight.
    To that end we construct a $3n-7$ dimensional set in $\ogr(2,n)\times \ogr(2,n)$ such that the Jacobian $\operatorname{Jac}_{\Psi_n}$ has full rank.
    Consider the following family of $2\times n$ matrices $A$ and $B$
    \begin{align}
        A = \begin{bmatrix}
            \nu_1  & \dots & \nu_n\\
           \nu_1 x_1 &  \dots & \nu_nx_n 
        \end{bmatrix},
        \quad
        B = \begin{bmatrix}
            \mu_1 &  \dots & \mu_n\\
            \mu_1y_1  & \dots & \mu_ny_n 
        \end{bmatrix}.
    \end{align}
    Then, by the $\operatorname{PGL}_2$ action on $A$ and $B$, we can assume that $x_1=0=y_1$, $x_2=1=y_2$ and $x_3=2=y_3$.
    Fixing in addition $\nu_i=1$ for all $i$ and $\mu_1=1$, we get a $(3n-7)$-dimensional family $\mathcal{V}\subset \gr(2,n)\times \gr(2,n)$ of pairs $(B,A)$.
    Denote the $2\times 2$ matrices with columns $(1,x_i)^T$ and $(1,x_j)^T$ by $(x_i,x_j)$, similarly for $(y_i,y_j)$.
    Then, for all $1\leq i<j\leq n$ we get
    \begin{align}
        z_{ij}=\frac{\det(\lambda_i y_i,\lambda_j y_j)}{\det(x_i,x_j)}
        =\lambda_i\lambda_j\frac{y_j-y_i}{x_j-x_i},
    \end{align}
    where $\lambda_i=\mu_i\in\C^\times$ and $\lambda_1=1$.
    Hence, the free parameters are $\lambda_2,\ldots,\lambda_n$,
    $x_4,\ldots,x_n$ and $y_4,\ldots,y_n$, giving
    $(n-1)+2(n-3)=3n-7$ free parameters.
    We will use the following $3n-7$ cross-ratios on $\P^N$ defined for all $j\geq 4$:
    \begin{align}
        u_j &= \frac{z_{1j}z_{23}}{z_{12}z_{13}}, \quad u_2=\frac{z_{23}}{z_{13}} , \quad u_3 = \frac{z_{23}}{z_{12}},\\
        v_j &= \frac{z_{2j}z_{13}}{z_{1j}z_{23}},\\
        w_j &= \frac{z_{3j}z_{12}}{z_{1j}z_{23}}.
    \end{align}
    These cross ratios give rise to a morphism $\phi: \mathring{\pen_n}\rightarrow \C^{3n-7}_{(u,v,w)}$.
    Note with our choice of $\mathcal{V}$ from above, these cross ratios simplify on $\Psi_n(\mathcal{V})$ to
    \begin{align}
         u_j &= \lambda_j\frac{y_j}{x_j}, \quad u_2=\lambda_2 , \quad u_3 = \lambda_3,\\
        v_j &= \frac{x_j(y_j-1)}{y_j(x_j-1)},\\
        w_j &= \frac{x_j(y_j-2)}{y_j(x_j-2)}.
    \end{align}
    Moreover, note that each ratio only depends on $x_j$, $y_j$, and $\lambda_j$ for the same index $j$, so in particular
    \begin{align}\label{eq:j_dependence}
        &\frac{\partial u_j}{\partial x_i}=0 ,\quad \frac{\partial u_j}{\partial y_i}=0 ,\quad \frac{\partial u_j}{\partial \lambda_i}=0, \\
        &\frac{\partial v_j}{\partial x_i}=0 ,\quad \frac{\partial v_j}{\partial y_i}=0 ,\quad \frac{\partial v_j}{\partial \lambda_i}=0, \\
        &\frac{\partial w_j}{\partial x_i}=0 ,\quad \frac{\partial w_j}{\partial y_i}=0 ,\quad \frac{\partial w_j}{\partial \lambda_i}=0,
    \end{align}
    for all $i\neq j$.
    Next consider the Jacobian $J=\operatorname{Jac}_{\phi\circ\Psi_n}$ of the composition of the crossratio and the Hadamard ratio $\phi\circ\Psi_n:\C^{3n-7}_{(\lambda,x,y)}\dashrightarrow \C^{3n-7}_{(u,v,w)}$, and note $\rk J\leq \rk \operatorname{Jac}_{\Psi_n}$, hence it suffices to show that $\rk J= 3n-7$. 
    We compute the rank of $J$ with respect to an ordered basis $(\lambda_2,\lambda_3,\lambda_4,x_4,y_4,\ldots, \lambda_n, x_n,y_n)$ of $\C^{3n-7}_{(\lambda,x,y)}$ and $(u_2,u_3,u_4,v_4,w_4,\ldots, u_n,v_n,w_n)$ of $\C^{3n-7}_{(u,v,w)}$.
    Then, $J$ is block diagonal with
    \[
        J = 1 \oplus 1 \oplus \bigoplus_{i=4}^nJ_i,
    \]
    by \eqref{eq:j_dependence} where 
    \[
    J_i = \begin{bmatrix}
        \frac{\partial u_i}{\partial\lambda_i} & \frac{\partial u_i}{\partial x_i} &\frac{\partial u_i}{\partial y_i} \\
        \frac{\partial v_i}{\partial\lambda_i} & \frac{\partial v_i}{\partial x_i} &\frac{\partial v_i}{\partial y_i} \\
        \frac{\partial w_i}{\partial\lambda_i} & \frac{\partial w_i}{\partial x_i} &\frac{\partial w_i}{\partial y_i} 
    \end{bmatrix} \quad\text{for all } i\geq 4.
    \]
    Then, note that for all $i$
    \[
    \det J_i=2\frac{y_i-x_i}{y_i^2(x_i-1)^2(x_i-2)^2},
    \]
    hence
    \[
    \det J = \prod_{i=4}^n\det J_i=2^{n-3}\prod_{i=4}^n\frac{y_i-x_i}{y_i^2(x_i-1)^2(x_i-2)^2}.
    \]
    For a generic choice of $(\lambda,x,y)\in (\C^\times)^{n-1}\times \C^{n-3}\times \C^{n-3}$, this determinant will not vanish, hence $J$ has generically full rank $3n-7$.
    This shows that $\rk\operatorname{Jac}_{\Psi_n}\geq \rk J=3n-7$.
    Therefore, $\dim \im \Psi_n\geq \dim \im \Psi_n\vert_{\mathcal{V}}\geq 3n-7$.
    Since $\dim \Pen_n\leq 3n-7$ by \eqref{eq:upper_bound_dim}, the set $\Psi(\mathcal{V})$ is dense in $\pen_n$.
    Therefore, $\dim \pen_n=3n-7$ as had to be shown.
    %By the inverse function theorem, the image $\mathring\pen_n$ of the restriction of $\Psi_n$ contains a Zariski open subset in $\pen_n$ of dimension $3n-7$. This proves that $\dim \mathring\pen_n\geq 3n-7$, concluding the proof.
\end{proof}

\begin{remark}
    Proposition \ref{prop:dimPen} also shows that the generic fiber of $\Psi_n$ has dimension $n-1$. Note that there is an $(n-1)$-dimensional algebraic torus contained in each fiber as follows. On the Pl\"ucker coordinates on each of the two factors, there is an action of the $(n-1)$-dimensional torus $(\C^\times)^{n-1}$, rescaling both $p_{ij}$ and $q_{ij}$ by $\lambda_i\lambda_j$. Their ratio is constant and the fiber over each point contains $(\C^\times)^{n-1}$.
    Since the generic fiber is $(n-1)$-dimensional, an irreducible component thereof  is a Zariski closure of $(\C^\times)^{n-1}$. 
\end{remark}

\begin{comment}
Let us also record the following technical results for later use.

\begin{lemma}\label{lem: containment I_n inside P_n}
Let $P_n$ be the vanishing prime ideal of the irreducible and reduced variety $\Pen_n$. Then we have the containment $P_n\supseteq I_n$. 
\begin{proof}
The proof of Proposition \ref{prop:dimPen} shows that there exists an open subset $V'\subset V\cap \ogr(2,n)\times \ogr(2,n)$ where the restriction of $\Psi_n$ is a morphism, $\Psi_n(V')$ has dimension $3n-7$ and $\Psi_n(V')\subset \Pen_n\setminus \lbrace z_{ij}=0 \ |  \ 1\leq i<j\leq n\rbrace$.  
Thus the closure $\overline{\Psi_n(V')}\subset \overline{\Pen_n\setminus \bigcup_{1\leq i < j\leq n} \lbrace z_{ij}=0\rbrace}=\overline{\mathbb{V}(P_n)\setminus \bigcup_{1\leq i < j\leq n} \lbrace z_{ij}=0\rbrace}$. Since $\Pen_n=\mathbb{V}(P_n)\subset \mathbb{V}(J_n)$, one has that the 
\[
\overline{\Psi_n(V')}\subset 
\overline{\mathbb{V}(J_n)\setminus \bigcup_{1\leq i < j\leq n} \lbrace z_{ij}=0\rbrace} = V\biggl(J_n:\,\bigl[\prod_{1 \leq i < j \leq n} \!\!\!  z_{ij} \,\bigr]^\infty \,\biggr) = \mathbb{V}(I_n).
\]
Since $\Psi_n(V')$ has dimension $3n-7$ and it is irreducible, its closure is
$\overline{\Psi_n(V')}=\Pen_n$, because they are both irreducible and of the same dimension. Thus $P_n\supseteq I_n$.
\end{proof}
\end{lemma}

\begin{remark}
In the course of the proof of the previous lemma, we notice that $\overline{ \mathbb{V}(J_n)\cap T}=\mathbb{V}(I_n)$, where  $T=\{z_{ij}\neq 0 \text{ for all } i<j\}$ is the algebraic torus. 
\end{remark}
\end{comment}

The pentad variety and the MHV gravity amplitudes can be computationally connected via the following matrix-completion property of the pentad from \eqref{eq:pentad}.
\begin{proposition}\label{prop: lifting to rank 2 symm matrices n=5}
    For $n=5$, there exists a Zariski dense subset of the affine cone over the pentad variety $\pen_5$,  where each point lifts to a unique symmetric matrix $\mathcal{Z}_5$ with the following propety: its diagonal entries $z_{ii}$ are given by rational functions $f_i=g_i/h_i$, $g_i,h_i\in \C[z_{ij}]$, in the off-diagonal entries $z_{ij}$ in such a way that $\rk \mathcal{Z}_5\leq 2$.
\end{proposition}
\begin{proof}
\begin{comment}
    This can be directly computed via elimination in \texttt{Macaulay2}. 
    More precisely, the claim is that, if $P$ is the ideal generated by the Pl\"ucker relations in $R_5$, then
    \[
        \bigg(P\cap (p_{ij}z_{ij}-q_{ij}\mid 1\leq i<j\leq 5)\bigg)\cap \C[z_{ij}\mid 1\leq i<j\leq 5] = (\operatorname{pen}_{12345}),
    \]
    in graded reverse lexicographic order.
\end{comment}
Introduce the following rational functions in the off-diagonal entries: 
\begin{align}
    f_1 &= \frac{z_{12}z_{14}z_{35}-z_{12}z_{15}z_{34}-z_{13}z_{14}z_{25}+z_{13}z_{15}z_{24}}{z_{24}z_{35}-z_{25}z_{34}},\\
    f_2 &= \frac{z_{12}z_{23}z_{45}-z_{12}z_{24}z_{35}+z_{13}z_{24}z_{25}-z_{14}z_{23}z_{25}}{z_{13}z_{45}-z_{14}z_{35}},\\
    f_3 &= -\frac{z_{13}z_{24}z_{35}-z_{13}z_{25}z_{34}-z_{14}z_{23}z_{35}+z_{15}z_{23}z_{34}}{z_{14}z_{25}-z_{15}z_{24}},\\
    f_4 &= \frac{z_{12}z_{34}z_{45}-z_{13}z_{24}z_{45}+z_{14}z_{24}z_{35}-z_{14}z_{25}z_{34}}{z_{12}z_{35}-z_{13}z_{25}},\\
    f_5 &= \frac{z_{13}z_{25}z_{45}-z_{14}z_{25}z_{35}-z_{15}z_{23}z_{45}+z_{15}z_{24}z_{35}}{z_{13}z_{24}-z_{14}z_{23}}.
    \end{align}
Let 
\[
\mathcal{Z}_5 =
    \begin{bmatrix}
        f_1 & z_{12} & z_{13} & z_{14} & z_{15} \\
        z_{12} & f_2 & z_{23} & z_{24} & z_{25} \\
        z_{13} & z_{23} & f_3 & z_{34} & z_{35} \\
        z_{14} & z_{24} & z_{34} & f_4 & z_{45} \\
        z_{15} & z_{25} & z_{35} & z_{45} & f_5
    \end{bmatrix}.
    \]
Let $U$ be the Zariski dense subset in the affine space of symmetric $5\times 5$ matrices in the $z_{ij}$'s (including the diagonal entries), where all the $f_i$'s make regular functions. 
One can check in \texttt{Macaulay2} that each $3\times 3$ minor of $\mathcal{Z}_5$ is either identically zero or it is a multiple of the pentad $\text{pen}_{12345}$. Now, let project $U$ to the affine space of the off-diagonal entries and obtain an open set $V$. Letting $W=V\cap \widehat{\Pen_5}\neq \emptyset$, the intersection of $V$ with the affine cone over $\Pen_5$, we see that any point in $W$ is such that its lift $\mathcal{Z}_5\in U$ has
$\rk \mathcal{Z}_5=2$. By the construction of the $f_i$'s, such a lift is unique on $W$. 
\end{proof}

\begin{remark}
It is easy to see that the uniqueness statement fails on on the entire variety $\Pen_5$. 
\end{remark}

We end this subsection by recording a technical result that will be useful later to prove our main theorem.
\begin{lemma}
    Recall that $Z_n =\C[z_{ij}\mid 1\leq i<j\leq n]$ and let $Z_n^{\pm}$ its Laurent polynomial ring.
    Let $J^p_n$ and $I^p_n$ be the ideals in $Z_n$ defined in \eqref{eq:def_Jpn} and \eqref{eq:def_Ipn}.
    Then, $I^p_n=
    \bigl(J^p_n\cdot Z_n^{\pm}\bigr)\cap Z_n$.
    %, where $J^p_n\cdot Z^{\pm}_n$ is the localization of $J^p_n$ in $Z^{\pm}_n$.
    The same holds true for the radicals, i.e. $\sqrt{I^p_n}=
    \bigl(\sqrt{J^p_n}\cdot Z_n^{\pm}\bigr)\cap Z_n$.
\end{lemma}
\begin{proof}
    Let $m=\prod_{1\le i<j\le n} z_{ij}$. Define the multiplicative set $S=\lbrace m^i, i\geq 0\rbrace$. 
    Then the localization of $Z_n$ at $S$ is
    \[
    S^{-1}Z_n = Z_n[m^{-1}] = Z^{\pm}_n,
    \]
    which is the coordinate ring of the torus $T=\{z_{ij}\neq 0 \text{ for all } i<j\}$. Let $J^p_n\subseteq Z_n$ be the ideal generated by the pentads, and recall that
    \[
    I^p_n = \bigl(J^p_n:m^\infty\bigr)
    =
    \{\, f\in Z_n \mid \exists k\ge 0 \text{ such that } m^k f\in J^p_n \,\}.
    \]
    We claim that
    \[
    I^p_n
    =
    \bigl(J^p_n\cdot Z^{\pm}_n\bigr)\cap \mathbb{C}[z_{ij}].
    \]
    In other words, $I^p_n$ is exactly the contraction to $Z_n$ of the extension of $J^p_n$ to the torus coordinate ring. Let $f\in I^p_n$. By definition of saturation, there exists $k\ge 0$ such that $m^k f\in J^p_n$. 
    Passing to the localization $Z_n[m^{-1}]$, we may divide by $m^k$, and obtain
    \[
    f=\frac{m^k f}{m^k}\in J^p_n\cdot Z_n[m^{-1}]
    =J^p_n\cdot Z^{\pm}_n.
    \]
    Since already $f\in Z_n$, it follows that
    $f\in \bigl(J^p_n\cdot Z^{\pm}_n\bigr)\cap \mathbb{C}[z_{ij}]$. 

    Now we show the converse, so assume that $f\in \bigl(J^p_n\cdot Z^{\pm}_n\bigr)\cap \mathbb{C}[z_{ij}]$.  
    So in the localized ring $Z_n[m^{-1}]$ we can write
    $f=\sum_{\ell=1}^r g_\ell h_\ell, g_\ell\in J^p_n,\ h_\ell\in Z_n[m^{-1}]$. 
    Each $h_\ell$ has the form $h_\ell=\frac{a_\ell}{m^{N_\ell}}$ where $a_\ell\in Z_n$. 
    Choose $N$ large enough so that $N\ge N_\ell$ for all $\ell$. Then
    \[
    m^N f
    =
    \sum_{\ell=1}^r g_\ell\, a_\ell\, m^{N-N_\ell}.
    \]
    Since $g_\ell\in J^p_n$ and $J^p_n$ is an ideal of $Z_n$, every term on the right-hand side lies in $J^p_n$. Therefore $m^N f\in J^p_n$. 
    This implies that $f\in J^p_n:m^\infty=I^p_n$.
    The proof for the statement about the radical ideals is identical, after replacing $f$ with a suitable power $f^\ell$.
\end{proof}

\begin{remark}\label{rem:contraction_prime}
    Suppose that the localized ideal $J^p_n\cdot Z_n^{\pm}$
    is a prime ideal of $Z_n^{\pm}$. Then its contraction $I^p_n =
    \bigl(J^p_n\cdot Z_n^{\pm}\bigr)\cap Z_n$
    is prime in $Z_n$. Indeed, if $ab\in I^p_n$ for some $a,b\in Z_n$, then
    $ab\in J^p_n\cdot Z_n^{\pm}$. 
    Since the latter ideal is prime, either
    $a\in J^p_n\cdot Z^{\pm}_n$ or
    $b\in J^p_n\cdot Z^{\pm}_n$. 
    Because $a,b\in Z_n$, this implies $a\in I^p_n$ or $b\in I^p_n$. Hence $I^p_n$ is prime.
\end{remark}

\subsection{Affine Slope Variety}\label{sec:affine_slope}

The purpose of this section is to prove our main result, Theorem \ref{thm:main_pentad}, on the pentad variety.
In order to do so, we define affine slope varieties as first introduced by Martin in \cite{martin2006slopesdeterminednpoints}.
It turns out that there is a nice connection between these affine slope varieties and the pentad via a Hadamard product with a toric variety, which plays the role of the little-group scaling, see Corollary \ref{cor:main_aff_to_pen}.
Moreover, we also realize the affine slope variety as the Hadamard ratio of two copies of the moduli space of $n$ marked points on $\P^1$, see Proposition \ref{prop:Ratio_M0n_slope}.

\medskip

\begin{definition}[Affine slope variety]\label{def:affine_slope}
    Consider $\aletter\geq2$ distinct points $p_1,\ldots, p_\aletter$ with coordinates $(x_i,y_i)$ in the affine complex plane $\A^2$, such that none of the $\binom{h}{2}$ lines $L_{ij}$ connecting any pair of points $p_i$ and $p_j$ is vertical.
    Let $m_{ij}=\frac{y_j-y_i}{x_j-x_i}$ be the indeterminates given by the slopes of the lines $L_{ij}$.
    The Zariski closure of all slope vectors $(m_{12},\ldots,m_{\aletter-1\aletter})\in \A^N$ with $N=\binom{\aletter}{2}$ is the algebraic variety $\mathcal{S}_h\subset \A^N$ called \emph{the affine slope variety}. 
\end{definition}

There is an equivalent description of the affine slope variety, also due to \cite{martin2006slopesdeterminednpoints}, in terms of graphs, which will be useful to us.
Let $\vzero, \vone,\dots, \vh$ be distinct vertices.
The $(\aletter{-}1)$-\emph{wheel} $W=W(\vzero;\vone,\dots, \vh)$ is defined as the graph with the \emph{edge set}
\[
   E(W)\coloneqq \{\vzero\vone,\vzero\vtwo,\dots, \vzero\vh\}\cup\{\vone\vtwo,\vtwo\vthree\dots, \vh\vone\}.
\]
Next consider the $(\aletter{-}1)\times (\aletter{-}1)$ matrix
\[
    M_{\vzero,\vone,\dots, \vh} \coloneqq
    \begin{bmatrix}
        m_{\vzero\vone} - m_{\vone\vtwo} & m_{\vone\vtwo} - m_{\vzero\vtwo}   & 0 &\dots&0\\
        0 & m_{\vzero\vtwo} - m_{\vtwo\vthree} & m_{\vtwo\vthree} - m_{\vzero\vthree}&\dots&0\\
        \vdots & \vdots & \vdots &\ddots&\vdots\\ m_{\vone\vh} -m_{\vzero\vone} & 0 & 0&\dots&  m_{\vzero\vh}-m_{\vone\vh} 
    \end{bmatrix}.
\]
% \[
%     M_{W_h} =
%     \begin{bmatrix}
%         m_{12} - m_{23} & m_{23} - m_{13}   & 0 &\dots&0\\
%         0 & m_{13} - m_{34} & m_{34} - m_{14}&\dots&0\\
%         \vdots & \vdots & \vdots &\ddots&\vdots\\ m_{2h} -m_{12} & 0 & 0&\dots&  m_{1h}-m_{2h} 
%     \end{bmatrix},
% \]
We adopt the convention that $v_{\aletter+i}=v_{i+1}$ for all $1\leq i\leq h-1$.
We call the determinants 
\begin{align}\label{eq:tree_poly}
\tau_{\vzero,\vone,\dots, \vh}\coloneqq \det M_{\vzero,\vone,\dots, \vh} =
\prod_{i=2}^h(m_{\vzero,v_i}-m_{v_i,v_{i+1}})-\prod_{i=2}^h(m_{\vzero,v_{i+1}}-m_{v_i,v_{i+1}}) 
\end{align}
% \begin{align*}
%     \tau(W) \coloneqq \det M_{W_h} 
%     =
%     \prod_{i=2}^h\big(m_{1i}-m_{i(i+1)}\big)-\prod_{i=2}^h\big(m_{1(i+1)}-m_{i(i+1)}\big) 
% \end{align*}
the \emph{tree-polynomial of the wheel $W=W(\vzero;\vone,\dots, \vh)$}, and sometimes also write $\tau(W)$ for $\tau_{\vzero,\vone,\dots, \vh}$.

\begin{example}\label{exm:tree_K4}
    Fix a set of vertices labeled as $\{1,2,3,4\}$.
    Then, the $3$-wheel $W_3=W(1;2,3,4)$ is the complete graph $K_4$.
    The matrix from above reads
    \[
        M_{W_3} =
    \begin{bmatrix}
        m_{12} - m_{23} & m_{23} - m_{13}   & 0 \\
        0               & m_{13} - m_{34}   & m_{34} - m_{14}\\
        m_{24} - m_{12} & 0                 & m_{14} - m_{24} 
    \end{bmatrix}.
    \]
    %where we write $m_{ij}=m_{v_iv_j}$ for brevity.
    We obtain the tree polynomial of $W_3$ as
    \begin{align}
        \tau(W_3)=\det M_{W_3} 
        &= -m_{12}m_{13}m_{24}+m_{12}m_{13}m_{34}+m_{12}m_{14}m_{23}-m_{12}m_{14}m_{34} \\  
        &\quad -m_{12}m_{23}m_{34}+m_{12}m_{24}m_{34}-m_{13}m_{14}m_{23}+m_{13}m_{14}m_{24}\\
        &\quad+m_{13}m_{23}m_{24}-m_{13}m_{24}m_{34}-m_{14}m_{23}m_{24}+m_{14}m_{23}m_{34}.
    \end{align}
    Note that this polynomial is square-free and has $12$ monomials much like the pentad constraint in \eqref{eq:pentad}.
\end{example}

We fix some more notation.
Denote by $A_h$ the tensor product
\[
    A_h\coloneqq \C[\lambda_i\mid 1\leq i\leq \aletter] \otimes_\C \C[m_{ij} \mid 1\leq i<j\leq \aletter].
\]
Recall that $Z_n=\C[z_{ij}\mid 1\leq i<j\leq n]$.
Denote by $A^{\pm}_\aletter$ the Laurent polynomial ring in the variables of $A_\aletter$, analogously for $Z_n^{\pm}$.
We then define the following two ideals in $\C[m_{ij}\mid i<j\leq \aletter]$.
Let $J_\aletter^s$ be the ideal generated by the tree polynomial of any $3$-wheel subgraph of the complete graph $K_\aletter$ on $\aletter$ vertices, i.e.,
\begin{align}\label{def:J^s_n}
    J_\aletter^s 
    &\coloneqq  \langle \tau(W_3) \mid W_3 \text{ wheel subgraph of } K_\aletter\rangle \\
    &\:= \langle \det M_{ijkl} \mid  W_3 = W(i;j,k,l)\text{ with }\{i,j,k,l\}\subset [\aletter]\,\rangle \subset \C[m_{ij}\mid 1\leq i<j\leq \aletter].
\end{align}
Similarly, let $I_\aletter^s$ be the ideal in $\C[m_{ij}\mid 1\leq i<j\leq \aletter]$ generated by all tree polynomials of any wheel subgraph of $K_\aletter$.

\begin{example}\label{exm:slope_ideals}
    For $h=5$ the affine slope variety $\ms_5\subset \C^{10}$ is of dimension $7$, the ideal $J_5^s$ has $5$ generators, which correspond to the $5$ complete graphs $K_4$ in $K_5$.
    The ideal $I_5^s$ in turn has $15$ extra generators corresponding to the $15=\frac{(h-1)!}{2}\binom{h}{4}$ wheels $W_4$ of size $4$ contained in $K_5$.
    However, the five generators of $J_5^s$ minimally generate $I_5^s$.
\end{example}

The following is the main result \cite[Theorem 1.1]{martin2006slopesdeterminednpoints} and a remark from the end of \emph{cit. loc}.

\begin{theorem}[{\cite{martin2006slopesdeterminednpoints}}]\label{thm:affine_slope_var}
    The affine slope variety $\mathcal{S}_h\subset \A^{\binom{\aletter}{2}}$ is irreducible of affine dimension $2\aletter-3$.
    Its prime vanishing ideal in $\C[m_{ij}\mid i<j\leq h]$ is $\slopeI{\aletter}$.
    Moreover, for $h\leq 9$ its vanishing ideal is equal to the ideal generated by the tree polynomials of all $3$-wheels $W_3$ contained in $K_\aletter$ i.e.
    \[
        \slopeI{\aletter} = \slopeJ{\aletter} \subset \C[m_{ij}\mid 1\leq i<j\leq h].
    \]
\end{theorem}

\begin{remark}
    Note that the affine slope variety $\mathcal{S}_\aletter\subset \C^{\binom{\aletter}{2}}$ is a cone.
    More precisely, for $t\in \C^\times$, and $m\in \mathcal{S}_\aletter$, we have $tm\in \mathcal{S}_\aletter$ as $t\bullet (x_i,y_i)\mapsto(x_i,t y_i)$ scales all $m_{ij}$ by $t$.
    Therefore, let $\P(\mathcal{S}_h)\subset \P^{\binom{h}{2}-1}$ be the \emph{projectivized slope variety}, which is the projectivization of the affine cone i.e. $(\ms_h\setminus \{0\})/\C^\times$.
    Equivalently the projectivization is simply given by taking the vanishing locus of $\slopeI{\aletter}$ in $\P^{\binom{h}{2}-1}=\P(\C^{\binom{h}{2}})$ as all generators are homogeneous.
    Conversely, taking the affine cone of the projectivization again recovers the affine slope variety.
    We will use these facts freely throughout the text.
\end{remark}

First let us connect affine slope varieties and Hadamard ratios.

\begin{proposition}\label{prop:Ratio_M0n_slope}
    Let $M_{0,n}$ be the moduli space of $n$-marked points on $\P^1$.
    Then, the Hadamard ratio $M_{0,n+1}:M_{0,n+1}$ is the projectivized slope variety $\P(\ms_{n})$.
\end{proposition}
\begin{proof}
    We construct a homogenized parametrization of $M_{0,n+1}$ as follows.
    Consider the $2\times (n{+}1)$ matrix 
    \begin{align}
        X =\begin{bmatrix}
            0 & 1 & \dots & 1 \\
            1 & x_1 & \dots & x_{n}
        \end{bmatrix},
    \end{align}
    in the indeterminates $x_i$ with $x_j\neq x_i$.
    This corresponds to a choice of $n$ marked affine points on $U_0\subset \P^1$ together with the point at infinity $[0:1]\in \P^1$.
    We fix the remaining affine coordinate freedom by setting $x_1=0$ and $x_n=1$.
    Taking the $2\times2$ minors of $X$ then embeds $M_{0,n+1}$ into $\P^M$, where $M=\binom{n+1}{2}-1$.
    Label the columns of $X$ from $0$ through $n$.
    For each $1\leq i < j \leq n$, denote by $p_{ij}=x_j-x_i$ the homogeneous $2\times 2$ minors of $X$ formed by columns $i$ and $j$, and note that $p_{0j}=-1$ for all $j>0$.
    Then, after projecting to the $\P^N$ with $N=\binom{n}{2}-1$ and coordinates $p_{ij}$ with $1\leq i < j\leq n$, the Hadamard ratio $M_{0,n+1}:M_{0,n+1}$ is the closure of the image of the map
    \[
        \Psi:M_{0,n+1}\times M_{0,n+1} \longrightarrow \P^N,\quad (q_{ij})\times (p_{ij}) \mapsto \frac{y_j-y_i}{x_j-x_i}=m_{ij}, \quad 1\leq i < j\leq n.
    \]
    After taking the affine cone over $\im \Psi$, we obtain a Zariski dense subset of the affine slope variety $\ms_{n}\subset \C^{N+1}$.
    Therefore, $M_{0,n+1}:M_{0,n+1}=\P(\ms_n)$, as the Hadamard ratio is defined by taking the Zariski closure of $\im \Psi$.
    This shows the claim.
    % Since $M_{0,n}$ is irreducible $\dim M_{0,n}=n-3$, the affine cone of the image of the map above has dimension $2n-5$, hence it is Zariski dense in $S_{n-1}$, which also has dimension $2n-5$.
\end{proof}

Recall that the tropicalization of the moduli space $M_{0,n}$ is precisely the space of all metric trees with $n$ labeled leaves.
Then, the following result addresses the relation between the affine slope variety and the BHV-phylogenetic tree space \cite{BILLERA2001733}, as first speculated in \cite{martin2006slopesdeterminednpoints}.

\begin{corollary}\label{cor:trop_aff_slope}
    The tropicalization of the projectivized slope variety $\P(\ms_{n})$ is given as the Minkowski difference of the tropicalization of two copies of $M_{0,n+1}$.
In symbols: 
    \[
        \trop(\P(\ms_{n})) = \trop(M_{0,n+1}) - \trop(M_{0,n+1}),
    \]
where $M_{0,n+1}$ is embedded as in Proposition \ref{prop:Ratio_M0n_slope}.
\end{corollary}
\begin{proof}
    This follows from combining Lemma \ref{lem:tropicalDiff} with Proposition \ref{prop:Ratio_M0n_slope}.
\end{proof}

The remainder of this section is devoted to prove Theorem \ref{thm:pentad_variety}.
We start by introducing a series of technical results.
To that end consider the ring isomorphism
\begin{align}\label{eq:laurent_iso}
    \phi: Z_n^{\pm 1}
    &\longrightarrow A_{n-1}^{\pm 1}\\
     z_{in}  &\longmapsto \lambda_{i}  \quad \quad \quad \:\: i<n,\\
     z_{ij} &\longmapsto \lambda_i\lambda_j m_{ij}  \quad i<j<n.
\end{align} 
Then, for the pentads $\operatorname{pen}_{ijkln}$, we have 
\begin{align}\label{eq:phi(pen)}
    \phi(\operatorname{pen}_{ijkln}) = - (\lambda_i\lambda_j\lambda_k\lambda_l)^2\det M_{ijkl}.
\end{align}

\begin{lemma}\label{lem:pen_as_pen_n}
    Let $n\geq 5$.
    The image of $\operatorname{pen}_{abcde}$ with $\{a,b,c,d,e\}\subset [n-1] $ under $\phi$ is contained in the following ideal
    \[
    \langle \phi(\operatorname{pen}_{ijkln}) \mid \{i,j,k,l\}\subset\{a,b,c,d,e\} \rangle
    \]
    in $\C[\lambda_i^{\pm1}\mid i\leq n-1]\otimes \C[m_{ij}^{\pm1}\mid i<j\leq n-1]$.
\end{lemma}
\begin{proof}

    It suffices to note we can write 
    \begin{align}
        \phi(\operatorname{pen}_{abcde}) 
        &=  \frac{1}{2}(\lambda_a\lambda_b\lambda_c\lambda_d\lambda_e)^2\bigg(A_{abcd}\det M_{abcd} + A_{abce}\det M_{abce}\\
        &+ A_{abde}\det M_{abde} + A_{acde}\det M_{acde} + A_{bcde}\det M_{bcde}\bigg),
    \end{align}
    with 
    \begin{align}
        A_{abcd}&=-2m_{ae}m_{be}
            +m_{ab}(m_{ae}+m_{be}-m_{ce}-m_{de}),\\
        A_{abce}&=2m_{ad}m_{bd}
            +m_{ab}(-m_{ad}-m_{bd}+m_{cd}+m_{de}),\\
        A_{abde}&=-2m_{ac}m_{bc}
            +m_{ab}(m_{ac}+m_{bc}-m_{cd}-m_{ce}),\\
        A_{acde}&=m_{ab}(m_{bc}+m_{bd}+m_{be}-m_{ab}),\\
        A_{bcde}&=m_{ab}(m_{ab}-m_{ac}-m_{ad}-m_{ae}).
\end{align}
    Then the claim follows as $(\lambda_a\lambda_b\lambda_c\lambda_d\lambda_e)^2$ is a unit in $\C[\lambda_i^{\pm1}\mid i<n]$.

\end{proof}

This has an immediate consequence.
\begin{lemma}\label{lem:slope_ideal_pentad_ideal}
    Let $J_n^p$ be the ideal generated by all pentads on $[n]$.
    Then, the image of the localization $J^p_n\cdot\C[z_{ij}^{\pm1}\mid i<j\leq n\:]$ under $\phi$, is the localization of the ideal $J^s_{n-1}$ generated by all tree polynomials of $3$-wheels inside $K_{n-1}$, i.e.
    \begin{align}
        \phi(J_n^p\cdot Z_n^\pm) = \slopeJ{n-1}\cdot A^{\pm}_{n-1}.
    \end{align}
    The same claim holds for the radicals of the ideals, i.e. 
    \begin{align}
        \phi\left(\sqrt{J_n^p}\cdot Z_n^\pm\right) = \sqrt{\slopeJ{n-1}}\cdot A^{\pm}_{n-1}.
    \end{align}
\end{lemma}
\begin{proof}
    This follows from combining the description of the images of the pentad equations in \eqref{eq:phi(pen)} with Lemma \ref{lem:pen_as_pen_n}.
    If $f\in \sqrt{J^p_{n}}$ then $f^m\in J^p_n$ for a suitable $m\in \Z_{\geq 1}$, then the radical statement follows analogously.
\end{proof}

% \begin{lemma}
%     Let $x_{ij}=m_{0i} - m_{ij}$ and identify $m_{ij}=m_{ji}$ for $i>j$.
%     Let $W$ be the wheel $W(0,1\ldots,h)$ 
%     Let $W' = W(0;r,\ldots, h,1)$ and $W'' = W(0;1,\ldots, r)$.
%     Then, we have 
%     \[
%         x_{1r}\tau(W), x_{r1}\tau(W) \in \langle \tau(W') , \tau(W'') \rangle.
%     \]
% \end{lemma}

\begin{lemma}\label{lem:split_wheel}
    Let $x_{ij}=m_{1i} - m_{ij}$ and identify $m_{ij}=m_{ji}$ for $i>j$.
    Let $W$ be the wheel $W(\vzero;\vone\ldots,\vh)$.
    For each pair $(v_a,v_b)\not\in E(W)$, there exist two wheels $\subwheelone_{ab}$ and $\subwheeltwo_{ab}$ such that their vertex sets are proper subsets of $\{\vzero,\vone,\ldots, \vh\}$, and
    \begin{align}
      x_{v_a v_b}\tau(W), x_{v_b v_a}\tau(W) \in \langle \tau(\subwheelone_{v_av_b}) , \tau(\subwheeltwo_{v_av_b}) \rangle.  \label{equation_split_wheel}
    \end{align}
    In particular, if $h\geq 7$, then \[\tau(W)^{\aletter}\in \sum_{(v_a,v_b)\not\in E(W)}\langle \tau(\subwheelone_{v_a,v_b}) , \tau(\subwheeltwo_{v_a,v_b}) \rangle.\]
        
\end{lemma}
\begin{proof} 
    We start by splitting the $h-1$ wheel $W(v_1,v_2,\ldots,v_h)$ as follows.
    Let $\subwheelone = W(\vzero;\vone,\ldots, v_r)$ and $\subwheeltwo = W(\vzero;v_r,\ldots, \vh,\vone)$ with $4\leq r\leq h-1$ (see Figure~\ref{fig:wheel_split}).
    
     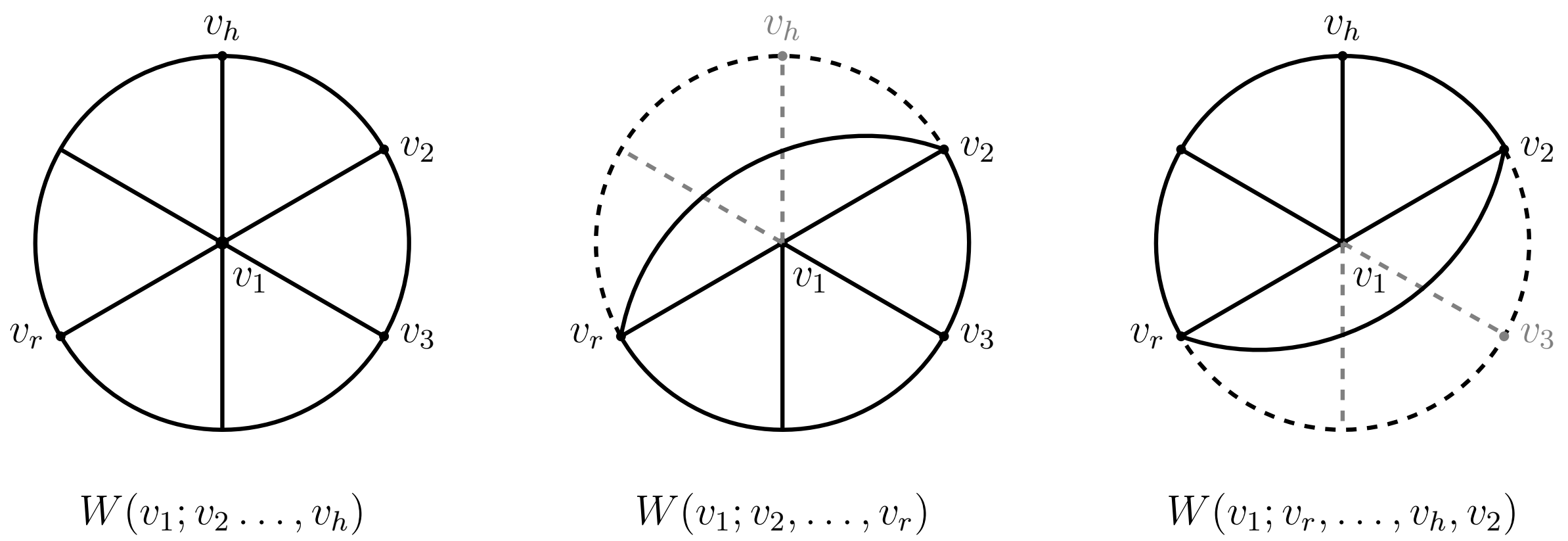
\begin{figure}[H]
         \begin{tikzpicture}[scale=0.47, very thick]

        %================== First graph ==================
        % Outer cycle
        \draw (3.47,-2) arc (-30:-150:4cm);
        \draw (-3.47,-2) arc (210:90:4cm);
        \draw[] (0,4) arc (90:-30:4cm);
        
        \node[below] at (0.6,-0.2) {$v_1$};
        \fill (0,0) circle (4pt);
        
        % Outer vertices
        \fill (90:4cm) circle (3pt);
        \node[above] at (90:4cm) {$v_\aletter$};
        \fill (30:4cm) circle (3pt);
        \node[right] at (30:4cm) {$v_2$};
        \fill (-30:4cm) circle (3pt);
        \node[right] at (-30:4cm) {$v_3$};
        %\fill (-90:4cm) circle (3pt);
        %\node[below] at (-90:4cm) {$\cdots$};
        \fill (-150:4cm) circle (3pt);
        \node[left] at (-150:4cm) {$v_r$};
        %\fill (150:4cm) circle (3pt);
        %\node[left] at (150:4cm) {$\iddots$};
        
        % Spokes
        \draw (0,0) -- (90:4cm);
        \draw (0,0) -- (30:4cm);
        \draw (0,0) -- (-30:4cm);
        \draw (0,0) -- (-90:4cm);
        \draw (0,0) -- (-150:4cm);
        \draw (0,0) -- (150:4cm);
        \node at (0,-5.8) {$W(\vzero;\vone\ldots,\vh)$};

        %================== Second graph ==================
        % Center
        \fill (12,0) circle (3pt);
        
        % Outer cycle
        \draw (15.46,2) arc (30:-150:4cm);
        \draw[dashed] (8.54,-2) arc (210:30:4cm);
        
        \fill (12,0) circle (3pt);      \node[below] at (12.6,-0.2) {$v_1$};
        \fill (15.46,2) circle (3pt);   \node[right] at (15.46,2) {$v_2$};
        \fill (15.46,-2) circle (3pt);  \node[right] at (15.46,-2) {$v_3$};
        %\fill (12,-4) circle (3pt);     \node[below] at (12,-4) {$\cdots$};
        \fill (8.54,-2) circle (3pt);   \node[left] at (8.54,-2) {$v_r$};
        \fill[gray] (12,4) circle (3pt);      \node[above] at (12,4) {\textcolor{gray}{$v_\aletter$}};
        %\fill (8.54,2) circle (3pt);   \node[left] at (8.54,2) {$\iddots$};
        \node at (12,-5.8) {$W(\vzero;\vone,\ldots, v_r)$};
        \draw (12,0)--(15.46,2);
        \draw (12,0)--(15.46,-2);
        \draw (12,0)--(12,-4);
        \draw[] (12,0)--(8.54,-2);
        \draw[dashed, gray] (12,0)--(8.54,2);
        \draw[dashed, gray] (12,0)--(12,4);
        \draw[bend right=50] (15.46,2) to (8.54,-2);

        %================== Third graph ===================
        \fill (24,0) circle (3pt);
        % Outer cycle
        \draw[dashed] (27.46,2) arc (30:-150:4cm);
        \draw[] (20.54,-2) arc (210:90:4cm);
        \draw[] (24.04,4) arc (90:30:4cm);
        \fill[gray] (27.46,-2) circle (3pt);  \node[right] at (27.46,-2) {\textcolor{gray}{$v_3$}};
        %\fill[gray] (24,-4) circle (3pt);     \node[below] at (24,-4) {};
        \fill (24,0) circle (3pt);      \node[below] at (24.6,-0.2) {$v_1$};
        \fill (24,4) circle (3pt);      \node[above] at (24,4) {$v_\aletter$};
        \fill (27.46,2) circle (3pt);   \node[right] at (27.46,2) {$v_2$};
        \fill (20.54,-2) circle (3pt);  \node[left] at (20.54,-2) {$v_r$};
        \fill (20.54,2) circle (3pt);   \node[left] at (20.54,2) {};
        \node at (24,-5.8) {$W(\vzero;v_r,\ldots, \vh,\vone)$};
        \draw (24,0)--(27.46,2);
        \draw[] (24,0)--(20.54,-2);
        \draw[] (24,0)--(24,4);
        \draw[] (24,0)--(20.54,2);
        \draw[dashed, gray] (24,0)--(27.54,-2);
        \draw[dashed, gray] (24,0)--(24,-4);
        \draw[bend left=50](27.46,2)to (20.54,-2);
        
        \end{tikzpicture}
        \caption{Splitting of the $h-1$ wheel $W(v_1;v_2\ldots,v_h)$.}
        \label{fig:wheel_split}
    \end{figure}

    \noindent First we show that
        \[
            x_{\vone v_r}\tau(W), x_{v_r\vone}\tau(W) \in \langle \tau(\subwheelone) , \tau(\subwheeltwo) \rangle.
        \]
    It suffices to consider the case $v_i=i$. Then, by definition, 
    \begin{align*}
    \tau(W)\!=\!\prod_{i=2}^h(m_{1i}-m_{i(i+1)})\!-\!\prod_{i=2}^h(m_{1(i+1)}\!-\!m_{i(i+1)})
    \!=\!\prod_{i=2}^h x_{i(i+1)}\!-\!\prod_{i=2}^h x_{(i+1)i}\!=\!AB\!-\!A'B',
    \end{align*}
    where $A=\prod_{i=2}^{r-1} x_{i(i+1)}$, $B=\prod_{i=r}^{h} x_{i(i+1)}$, $A'=\prod_{i=2}^{r-1} x_{(i+1)i}$, and $B'=\prod_{i=r}^h x_{(i+1)i}$.
    Moreover, \[\tau(\subwheelone)=\prod_{i=2}^r x_{i(i+1)}-\prod_{i=2}^r x_{(i+1)i}=Ax_{r2}-A'x_{2r},
    \]
    \[\tau(\subwheeltwo)=\Big(\prod_{i=r}^h x_{i(i+1)}\Big)x_{2r}-\Big(\prod_{i=r}^h x_{(i+1)i}\Big)x_{r2}=Bx_{2r}-B'x_{r2}.\]
    It follows that
    \begin{align*}    
        A\tau(\subwheeltwo)+B'\tau(\subwheelone)
        =ABx_{2r}-AB'x_{r2}+AB'x_{r2}-A'B'x_{2r}&\\
        =(AB-A'B')x_{2r}&=x_{2r}\tau(W).
    \end{align*}
    Similarly,
    \begin{align*}    
        A'\tau(\subwheeltwo)+B\tau(\subwheelone)
        &=A'Bx_{2r}-A'B'x_{r2}+ABx_{r2}-A'Bx_{2r}\\
        &=(AB-A'B')x_{r2}=x_{r2}\tau(W).
    \end{align*}
    Moreover, if $v_a$ and $v_b$ are non-adjacent in $W$ with $a<b$, then neither of them can be the center of the wheel since the center is adjacent to every other vertex. 
    Then \[W(\vzero;\vone, \dots, \vh)=W(\vzero;v_a, \dots, v_b,\dots, \vh,\vone,\dots, v_{a-1})=W(\vzero;\vone', \dots, \vh'),\] where $v_i'=v_{i+a-2}$.\\
    Now assume $W=W(1;2,\dots, h)$ is such that $h\ge 7$. 
    The containment relations in \eqref{equation_split_wheel} imply that
    $$ m_{1a}\tau(W)-m_{1b}\tau(W) =x_{ab}\tau(W)-x_{ba}\tau(W) \in \langle \tau(\subwheelone_{ab}) , \tau(\subwheeltwo_{ab}) \rangle, $$
    for each pair of non-adjacent vertices $a$ and $b$ of $W$. In particular,
    \[\begin{array}{@{}l@{\qquad}l@{}}
         (m_{12}-m_{1i})\tau(W)\in \langle \tau(\subwheelone_{2i}) , \tau(\subwheeltwo_{2i}) \rangle &\text{for all $4\le i\le h-1$},\\
         (m_{13}-m_{1j})\tau(W)\in \langle \tau(\subwheelone_{3j}) , \tau(\subwheeltwo_{3j}) \rangle &\text{for all $5\le j\le h$},\\
         (m_{14}-m_{1k})\tau(W)\in \langle \tau(\subwheelone_{4k}) , \tau(\subwheeltwo_{4k}) \rangle &\text{for all $k=6,\dots, h,2$},\\
         (m_{15}-m_{1l})\tau(W)\in \langle \tau(\subwheelone_{5l}) , \tau(\subwheeltwo_{5l}) \rangle &\text{for all $l=7,\dots, h, 2,3$},\\
         & \vdots\\
         (m_{1h}-m_{1\anotherletter})\tau(W)\in \langle \tau(\subwheelone_{\anotherletter h}) , \tau(\subwheeltwo_{\anotherletter h}) \rangle &\text{for all $3\leq \anotherletter\leq h-2$}.
    \end{array}\]
    It follows that $m_{1i}\tau(W)=m_{1j}\tau(W)$ in
    $\C[m_{ij}\vert 1\leq i<j\leq h]/H$ for all $2\le i<j\le h$, where
    \[
        H=\sum_{(v_a,v_b)\notin E(W)}\langle \tau(\subwheelone_{ab}) , \tau(\subwheeltwo_{ab}) \rangle\subseteq \C[m_{ij}\vert 1\leq i<j\leq h].
    \]
    Hence 
    \begin{align*}
        \tau(W)^{\aletter}
        =&\tau(W)^{h-1}\Big(\prod_{i=2}^h(m_{1i}-m_{i(i+1)})-\prod_{i=2}^h(m_{1(i+1)}-m_{i(i+1)})\Big)\\
        =&\prod_{i=2}^h(m_{1i}\tau(W)-m_{i(i+1)}\tau(W))-\prod_{i=2}^h(m_{1(i+1)}\tau(W)-m_{i(i+1)}\tau(W))\\
        =&\prod_{i=2}^h(m_{12}\tau(W)-m_{i(i+1)}\tau(W))-\prod_{i=2}^h(m_{12}\tau(W)-m_{i(i+1)}\tau(W))=0
    \end{align*}
    in $\C[m_{ij}\vert 1\leq i<j\leq h]/H$. That is, $\tau(W)^h\in H$. 
    Up to relabeling and cyclic rotation of the indices, it follows that, for every wheel $W=W(\vzero;\vone\ldots,\vh)$ with $h\ge7$,
    \[\tau(W)^h\in \sum_{(v_a,v_b)\not\in E(W)}\langle \tau(\subwheelone_{v_av_b}) , \tau(\subwheeltwo_{v_av_b}) \rangle.\qedhere\]
\end{proof}

For a homogeneous ideal $I$ in a standard graded polynomial ring, we write $I_{\leq d}$ for the ideal generated by all homogeneous polynomials in $I$ whose degrees are less than or equal to $d$. Then Lemma~\ref{lem:split_wheel} immediately implies the following.
\begin{corollary}\label{cor: reduce_degree_wheel}
    Let $h\geq 7$. Let $\tau(W)$ be a tree polynomial of an $i$-wheel $W$, which is a subgraph of $K_h$, then $\tau(W)^{i+1}\in \big(\slopeI{h}\big)_{\leq i-1}$ for all $6\leq i\leq h-1$. 
\end{corollary}
\begin{proof}
    By Lemma~\ref{lem:split_wheel},  
    \begin{align*}
        \tau(W)^{i+1}&\in \sum_{(v_a,v_b)\not\in E(W)}\langle \tau(\subwheelone_{v_a,v_b}) , \tau(\subwheeltwo_{v_a,v_b})\rangle \\
        &\subseteq \langle \text{tree polynomials of $j$-wheel subgraphs of $K_h$}\mid j\leq i-1\rangle  = \big(\slopeI{h}\big)_{\leq i-1}. \qedhere
    \end{align*} 
\end{proof}

The following lemma was proved in Enkosky's PhD thesis \cite[Theorem 4.1.2]{enkosky2011slope}. We provide here an alternative and simpler proof.
\begin{lemma}\label{lem:radical_equal}
    We have $\sqrt{J_h^s}= I_h^s$ for any $h\geq 4$.
\end{lemma}
\begin{proof}
    It is clear that $\slopeJ{h}\subseteq\slopeI{h}$ for any $h\geq4$, and equality holds for $h\leq 9$ by Theorem~\ref{thm:affine_slope_var}. Since $\slopeI{h}$ is a prime ideal, it follows that $\sqrt{\slopeJ{h}}\subseteq\sqrt{\slopeI{h}}=\slopeI{h}$ for any $h\geq 4$. It remains to prove the reverse inclusion for $h\geq 10$.
    
    Fix $h\geq 10$.
    We prove by induction on $i$ that ${\big(\slopeI{h}}\big)_{\leq i}\subseteq\sqrt{\slopeJ{h}}$ for all $3\leq i\leq h-1$.
    Since ${\big(\slopeI{h}}\big)_{\leq h-1}=\slopeI{h}$, this will complete the proof.
    For $3\leq i\leq 8$, the claim follows by applying Theorem~\ref{thm:affine_slope_var} again, which gives ${\big(\slopeI{h}}\big)_{\leq i}=\slopeJ{h}$. 
    Now let $9\leq i\leq \aletter-1$ and for the purpose of induction assume that ${\big(\slopeI{h}}\big)_{\leq i-1}\subseteq\sqrt{\slopeJ{h}}$.
    Let $W$ be an $i$-wheel subgraph of $K_h$.  
    By Corollary~\ref{cor: reduce_degree_wheel}, we have $\tau(W)^{i+1}\in {\big(\slopeI{h}}\big)_{\leq i-1}$.
    The inductive hypothesis therefore implies that
    $\tau(W)^{i+1}\in \sqrt{\slopeJ{h}} $, 
    which in turn implies $\tau(W)\in \sqrt{\slopeJ{h}} $.
    Hence every degree-$i$ generator  of $\slopeI{h}$ belongs to $\sqrt{\slopeJ{h}}$, and consequently ${\big(\slopeI{h}}\big)_{\leq i}\subseteq\sqrt{\slopeJ{h}}$. 
\end{proof}

We can finally prove the main result in this section.
\begin{proof}[Proof of Theorem 
\ref{thm:pentad_variety}]
    The dimensional claim is Proposition \ref{prop:dimPen}.
    Next, we show that $\sqrt{I^p_{n}}$ is prime.
    By Theorem~\ref{thm:affine_slope_var}, $\slopeI{n-1}$ is a prime ideal of $\C[m_{ij}\mid 1\leq i<j\leq n-1]$.
    By Lemma~\ref{lem:radical_equal}, it follows that $\sqrt{\slopeJ{n-1}}$ is also a prime ideal of $\C[m_{ij}\mid 1\leq i<j\leq n-1]$. We can extend this ideal to $A_{n-1}$ and obtain the prime ideal $\sqrt{\slopeJ{n-1}}\cdot A_{n-1}$.
    Now, we check that the monomial $u\coloneqq \prod_{1\leq i<j\leq n-1}m_{ij}\prod_{k=1}^{n-1} \lambda_k$ does not belong to $\sqrt{\slopeJ{n-1}}\cdot A_{n-1}$. 
Since $\sqrt{\slopeJ{n-1}}\cdot A_{n-1}$ is prime,
if $u\in \sqrt{\slopeJ{n-1}}\cdot A_{n-1}$, then 
at least one of the variables of $A_{n-1}$ must belong to $\sqrt{\slopeJ{n-1}}\cdot A_{n-1}$. However, by definition $\sqrt{\slopeJ{n-1}}=\slopeI{n-1}$ is generated by homogeneous polynomials of degree at least $3$. Consequently, its extension 
$\sqrt{\slopeJ{n-1}}\cdot A_{n-1}$ is also generated by homogeneous polynomials of degree at least $3$. In particular, it contains no nonzero homogeneous elements of degree $1$. This yields a contradiction.
\begin{comment}
Indeed, choose $x_i=y_i=i$ in the parameterization $m_{ij}=\frac{y_j-y_i}{x_j-x_i}$ for each $1\leq i<j\leq n-1$ and so we obtain $m_{ij}=1$;  choose $\lambda_k=1$ for all $1\leq k\leq n-1$.  With these choices, we have $u=1$. The point with coordinates $m_{ij}=1$ and $\lambda_k=1$ is then in the variety defined by the prime ideal $\sqrt{\slopeJ{n-1}}\cdot A_{n-1}$; so any polynomial in this ideal must vanish on this given point. However, $u=1$ at this point and therefore $u\notin \sqrt{\slopeJ{n-1}}\cdot A_{n-1}$. 
\end{comment}    
Since $u$ is not in the prime ideal $\sqrt{\slopeJ{n-1}}\cdot A_{n-1}$, localizing at $u$ preserves the primality of $\sqrt{\slopeJ{n-1}}\cdot A_{n-1}$. Hence 
$\sqrt{\slopeJ{n-1}}\cdot A^{\pm}_{n-1}$ is prime. Since $\phi$ is an isomorphism we have that 
\[
\sqrt{J^p_n}\cdot Z_{n}^\pm = \phi^{-1} \Big(\sqrt{\slopeJ{n-1}}\cdot A^{\pm}_{n-1}\Big)
\]
is a prime ideal of $Z_{n}^\pm$.
Moreover, Remark~\ref{rem:contraction_prime} shows that the primality of $\sqrt{J^p_n}\cdot Z_{n}^\pm $ forces $\sqrt{I^p_n}$ to be prime. Thus, we obtain the following chain of implications establishing primality:
    \[
        \slopeI{n-1} \text{ is prime }
        \Longrightarrow
        \sqrt{\slopeJ{n-1}} \text{ is prime }
        \Longrightarrow
        \sqrt{J^p_n}\cdot Z_{n}^\pm \text{ is prime }
        \Longrightarrow
        \sqrt{I^p_n} \text{ is prime }.
    \]
    In particular, if $n\leq 10$, then the same chain of implications shows that $I^p_n$ is prime:
    \[
         \slopeI{n-1} \text{ is prime }
         \Longrightarrow
         \slopeJ{n-1} \text{ is prime }
         \Longrightarrow
         J^p_n\cdot Z_{n}^\pm \text{ is prime }
         \Longrightarrow
         I^p_n \text{ is prime.}
     \]
     This concludes the proof.
    % By Theorem \ref{thm:affine_slope_var} the vanishing ideal of $I(\mathcal{S}_{n-1})$ is prime and has affine dimension $2n-5$ in $\C[m_{ij}\mid i<j\leq n-1]$.
    % Hence, the localization of $I(\mathcal{S}_{n-1})$ in $\C[\lambda_{i}^{\pm1}\mid i \leq n-1] \otimes \C[m_{ij}^{\pm1}\mid i<j \leq n-1]$ is prime of affine dimension $3n-6 = 2n-5 + n-1$, which matches the dimension of $J_n\C[z_{ij}^{\pm1}\mid i<j\leq n]$.
    % Next, recall from Theorem \ref{thm:affine_slope_var} that for $n\leq 10$ the ideal $I(S_{n-1})\C[m_{ij}^{\pm1}\mid i<j\leq n-1]$ is generated by the tree polynomials $\det M_{ijkl}$ for all $\{i,j,k,l\}\subset [n-1]$.
    % Then the claim follows from Lemma \ref{lem:slope_ideal_pentad_ideal} and Remark \ref{rem:contraction_prime}.
\end{proof}

We end this section by recording our main structural result relating the affine slope variety with the pentad variety, i.e. Corollary \ref{cor:main_aff_to_pen}.
We restate this result here as follows.

\begin{proposition}\label{prop:toric_times_slope_pentad}
    Let $\Delta$ be the second hypersimplex in $\R^n$, i.e. the convex hull of all vectors $e_i+e_j\in \R^{n}$ for $i<j\leq n$, where $e_i$ are the standard basis vectors of $\R^n$.
    Let $X_{\Delta}\subset \P^N$ with $N=\binom{n}{2}-1$ denote the embedded toric variety associated to $\Delta$.
    Then, the Hadamard product $\P(\mathcal S_{n})*X_{\Delta}$ equals to the pentad variety $\pen_n$ as subvarieties of $\P^N$.
\end{proposition}
\begin{proof}
    %Let $\mathring{\pen_n}$ be the open subvariety of $\pen_n\subset \P^N$ with all coordinate hyperplanes $z_{ij}=0$ removed.
    First recall that $X_{\Delta}\subset \P^N$ is defined as the Zariski closure of the monomial map
    \[
        (\C^\times)^{n}\rightarrow \P^{N}, \quad \lambda \mapsto (\lambda_i\lambda_j)_{i<j}.
    \]
    Next, recall that the Hadamard product $\P(\ms_{n})*X_{\Delta}$ is the closure of the image of the following map
    \[
        \P(\ms_{n}) \times X_{\Delta} \dashrightarrow \P^{N},\quad (m_{ij},\lambda_i\lambda_j)\mapsto (\lambda_i\lambda_jm_{ij})=\left(\lambda_i\lambda_j\frac{y_j-y_i}{x_j-x_i}\right).
    \]
    Lastly, recall the parametrization in the proof of Proposition \ref{prop:dimPen} via
    \[
        (\C^\times)^{n} \times \C^{n} \times \C^{n} \dashrightarrow \Pen_n \subset \P^N,\quad (\lambda,x,y) \mapsto \left(\lambda_i\lambda_j\frac{y_j-y_i}{x_j-x_i}\right).
    \]
    % Then, clearly $U\subset \Pen_n$ is also open and dense in $\P(\ms_n)*X_{\Delta}$.
    % Therefore, after taking closures we find $\pen_n\subset \P(\ms_n)*X_{\Delta}$, since both $\pen_n$ and $\P(\ms_n)*X_{\Delta}$ are closed by definition.
    % Conversely, the parametrization above also gives a dense open subset $U'$ in $\P(\ms_n)*X_{\Delta}$, where we set $m_{ij}=\frac{y_j-y_i}{x_j-x_i}$ as usual.
    % This set $U'$ is also contained in $\pen_n$, therefore, again after taking closures we find $\P(\ms_n)*X_{\Delta}\subset \pen_n$.
    Hence, the closures of the images of the two parametrizations, which define the varieties $\Pen_n$ and $\P(\ms_n)*X_{\Delta}$, agree.
    This shows the claim.
    % Conversely, for a generic affine slope vector $m\in \C^{N+1}$ and a generic point $\lambda\in (\C^\times)^n$, we can make two elements in the Grassmannians, represented by two $2\times n$ matrices $Y=(\lambda_i,\lambda_iy_i)^T$ and $X_{\Delta}=(1,x_i)^T$.
    % Their Pl\"ucker coordinates are $q_{ij}=\lambda_{i}\lambda_j(y_j-y_i)$ and $p_{ij}=x_j-x_i$.
    % Therefore, $q_{ij}/p_{ij}$ make points in $\pen_n$.
    % Hence, after taking closures we have $\P(\ms_n)*X_{\Delta}\subset \pen_n$, as $\lambda$ and $m$ were generic, hence the closure is $\P(\ms_n)*X_{\Delta}$.
    % We saw in the proof of Proposition \ref{prop:dimPen} that generically the Jacobian of this map has full rank, hence this map is a dense parametrization of $\mathring{\pen_n}$ by a dimension count.
    %Then, taking the closure of this map, i.e. the Hadamard product $\ms_{n}*X_{\Delta}$, gives the closure of $\mathring{\pen_n}$, i.e. $\pen_n$, which had to be shown.
\end{proof}

\subsection{Affine Slopes and Phirotopes}

In this section, we want to point out another surprising connection of the pentad and the affine slope variety with the combinatorial structure of phirotopes \cite{phirotopes}, which are a generalization of chirotopes.
We record here the more general definition.

\begin{definition}\label{def:phirotope}
    Let $E=\{1,\ldots,n\}$ be a finite index set. 
    A function
    \[
    \varphi \colon E^d \longrightarrow S^1 \cup \{0\}
    \]
    on all tuples of $E$ is called a \emph{rank-$d$ phirotope} if:
    
    \begin{enumerate}[itemsep = 0pt]
        \item[(a)] The function $\varphi$ is alternating.
    
        \item[(b)] For any two subsets $\tau=\{\tau_1,\ldots, \tau_{d-1}\}$ and $\mu=\{\mu_1,\ldots,\mu_{d+1}\} $
        of $E$, define
        \[
            \omega_k = (-1)^{k-1} \varphi(\tau,\mu_k) \varphi(\mu_1,\ldots,\widehat{\mu_k},\ldots,\mu_{d+1}),
        \]
        for $k=1,\ldots,d+1$. Then there exist
        $r_1,\ldots,r_{d+1}\in\mathbb{R}_{>0}$ such that
        \[
            \sum_{k=1}^{d+1} r_k\omega_k=0.
        \]
    \end{enumerate}
    We call $\varphi$ \emph{uniform} if its image on subsets of $E^d$ without repeated elements, is contained in $S^1$.
    We say $\varphi$ is \emph{realizable} if there exists a configuration of points $Z=(Z_1,\ldots,Z_n)\in \C^{d\times n}$ such that 
    \[
        \varphi(i_1,\ldots, i_d) = \omega(\det(Z_{i_1},\ldots, Z_{i_d})),
    \]
    where $\omega:\C\rightarrow S^1\cup \{0\}$ is the phase map with
    \[
        \omega(re^{ix}) = \begin{cases}
            0 &,r=0\\
            e^{ix} &,r\neq 0.
        \end{cases}
    \]
\end{definition}

\begin{example}
     Let us consider an explicit example of a rank-$2$ phirotope on $4$ elements, with $E=\{1,2,3,4\}$.
    We define the map $\varphi:E^2\to S^1\cup\{0\}$ by setting
    \begin{align}
        \varphi(1,2)&=1,
        &
        \varphi(1,3)&=i,
        &
        \varphi(1,4)&=\frac{2+i}{\sqrt5},\\
        \varphi(2,3)&=\frac{-1+i}{\sqrt2},
        &
        \varphi(2,4)&=\frac{1+i}{\sqrt2},
        &
        \varphi(3,4)&=1,
    \end{align}
    and set $\varphi(i,j)=-\varphi(j,i)$ for $i<j$, as well as $\varphi(i,i)=0$.
    Then, in particular $\varphi$ is alternating, so it satisfies condition (a) from Definition \ref{def:phirotope}.
    Next, we verify condition (b).
    For the four indices $1,2,3,4$, set
    \begin{align}
        \omega_1
        &=
        \varphi(1,2)\varphi(3,4)=1,\\
        \omega_2
        &=
        -\varphi(1,3)\varphi(2,4) =\frac{1-i}{\sqrt2},\\
        \omega_3
        &=
        \varphi(1,4)\varphi(2,3)=\frac{-3+i}{\sqrt{10}}.
    \end{align}
    Therefore, we can obtain positive solutions with $r_1=2$, $r_2=\sqrt{2}$, and $r_3=\sqrt{10}$, such that
    \begin{align}
        0=r_1\omega_1+r_2\omega_2+r_3\omega_3
        &=
        2+(1-i)+(-3+i).
    \end{align}
    Hence, $\varphi$ indeed defines a phirotope.
    Moreover, by construction $\varphi$ is uniform as $\im \varphi\subset S^1$.

     We now show that $\varphi$ is realizable.
     Consider the following four affine complex points
    \[
        z_1=0,\qquad
        z_2=1,\qquad
        z_3=i,\qquad
        z_4=2+i
    \]
    in the complex plane \(\mathbb C\simeq\mathbb R^2\).
    Using homogeneous coordinates, these correspond to the vectors
    \[
        Z_1=\begin{pmatrix}1\\0\end{pmatrix},
        \qquad
        Z_2=\begin{pmatrix}1\\1\end{pmatrix},
        \qquad
        Z_3=\begin{pmatrix}1\\i\end{pmatrix},
        \qquad
        Z_4=\begin{pmatrix}1\\2+i\end{pmatrix}
        \in\mathbb C^2.
    \]
    Then, indeed, the phases $\omega$ of the determinants of $Z_iZ_h$ are given as $\varphi(i,j)$, hence $\varphi$ is a uniform realizable phirotope of rank $2$.
\end{example}

\medskip

\begin{remark}
    The original reason for trying connect phirotopes with the pentad variety, was the observation made by Bernd Sturmfels and J\"urgen Richter-Gebert that the following realizability condition appearing in \cite[Theorem 5.2]{phirotopes}
    \[
        \sum_{\sigma\in \Sym(5)}\operatorname{sgn}(\sigma)\,
        \varphi(\sigma(1,2))^2
        \varphi(\sigma(2,3))^2
        \varphi(\sigma(3,4))^2
        \varphi(\sigma(4,5))^2
        \varphi(\sigma(5,1))^2
        = 0,
    \]
    precisely recovers the pentad constraint from \eqref{eq:pentad} after substituting $z_{ij}=(\varphi(i,j))^2$.
\end{remark}

% In case of realizable phirotopes of rank $d$ represented by an arrangement of points $P_1,\ldots, P_{n}$, we write $\varphi((P_i,P_j))$ for $\varphi(i,j)$.
Our main result in this paragraph makes the remark above rigorous.
More specifically, it relates the realizable phirotopes of rank $2$ on $n{+}1$ points to the affine slope variety $\ms_n$.

\begin{proposition}\label{prop:phiro_to_slopes}
    For a uniform rank $2$ phirotope $\varphi$ with index set $E=\{0,1,\ldots, n\}$,
    denote by $\widehat{\varphi}(i,j) \coloneqq \varphi(i,j)/(\varphi(0,i)\varphi(0,j)$ with $1\leq i<j\leq n$ its normalized phases.
    Let $\widehat{\varphi}^{\,2}\coloneqq (\widehat{\varphi}(1,2)^2,\ldots, \widehat{\varphi}(n-1,n)^2)\in \C^{\binom{n}{2}}$ be the vector of these coordinates.
    Define $\mathcal{R}_n\subset \C^{\binom{n}{2}}$ as the set of all such points obtained from uniform realizable phirotopes $\varphi$ on $E$, that is
    \[
        \mathcal{R}_n \coloneqq \left\{ \widehat{\varphi}^{\,2} \mid \varphi \text{ uniform realizable phirotope on } E \right\} \subset \C^{\binom{n}{2}}.
    \]
    Then, the complex Zariski closure of $\mathcal{R}_n$ coincides with the affine slope variety $\ms_{n}$.
\end{proposition}
\begin{proof}
    First note that by \cite[Lemma 2.3]{phirotopes} any realizable phirotope of rank $2$ on $n+1$ elements can be represented by a point configuration $(Z_0,\ldots,Z_n)$ in $\P^1$ of $\varphi$.
    By \emph{loc. cit.} we can also assume that $Z_0=(0,1)^T$, and $Z_i=b_i(1,a_i)^T$ for $i>0$, with $b_i\in\C^\times$ and pairwise distinct $a_i\in\C$.
    Thus, for the phases we get $\varphi(0,i)=-\omega(b_i)$ and $\varphi(i,j)=\omega(b_i)\omega(b_j)\omega(a_j-a_i)$.
    With $a_i=X_i+\sqrt{-1}Y_i$, where $X_i,Y_i\in \R$, the normalized squared phases satisfy
    \begin{align}
        \widehat{\varphi}(i,j)^2
        = \frac{a_j-a_i}{\overline{a_j-a_i}}=\frac{(X_j-X_i) + \sqrt{-1}(Y_j-Y_i)}{(X_j-X_i) - \sqrt{-1}(Y_j-Y_i)},
    \end{align}
    where $\overline{a_j-a_i}$ denotes the complex conjugate.
    This is well-defined, since the $a_i$ are pairwise distinct.
    Then, after a complex change of parameters setting $x_i\coloneqq X_i-\sqrt{-1}Y_i$ and $y_i\coloneqq X_i+\sqrt{-1}Y_i$, we get
    \[
        \widehat{\varphi}(i,j)^2=\frac{y_j-y_i}{x_j-x_i}.
    \]
    Therefore, each point $\widehat{\varphi}^{\,2}\in \mathcal{R}_n$ corresponds to a point in the affine slope variety $\ms_n$.
    Note that we have $x_i=\overline{y_i}$.
    It suffices to show that $\mathcal{R}_n$ is Zariski dense in $\ms_n$.
    To that end, let
    \[
        \{x_i=\overline{y_i}\text{ for all }i\}\subset \C^{2n},
        \quad \text{and }
        U\coloneqq \{x_i\neq x_j\text{ for all }i\neq j\}\subset \C^{2n}.
    \]
    Since $\{x_i=\overline{y_i} \text{ for all }i\}\cong\R^{2n}$ is complex Zariski dense in $\C^{2n}$, the intersection with $U$ is Zariski dense in $U$.
    The morphism
    \[
        U\longrightarrow \C^{\binom{n}{2}},\quad
        (x_i,y_i)\mapsto \left(\frac{y_j-y_i}{x_j-x_i}\right)_{i<j}
    \]
    has image with Zariski closure $\ms_n$, and its restriction to $\{x_i=\overline{y_i}\text{ for all }i\}\cap U$ has image $\mathcal{R}_n$.
    Hence, $\mathcal{R}_n$ is Zariski dense in $\ms_n$.
    This shows the claim.
\end{proof}

In fact we can see that the affine slope variety completely characterizes the realizable uniform phirotopes.

\begin{corollary}\label{cor:realizability}
    Let $\varphi$ be a uniform rank $2$ phirotope on $E=\{0,1,\ldots,n\}$.
    The following are equivalent:
    \begin{enumerate}
        \item[(i)] $\varphi$ is realizable.
        \item[(ii)] $\widehat{\varphi}^{\,2}\in\ms_n$.
        \item[(iii)] $\tau_{ijkl}(\widehat{\varphi}^{\,2})=0$ for every $1\leq i<j<k<l\leq n$.
        %\item[(iv)] The restriction of $\varphi$ to $\{0,i,j,k,l\}$ is realizable for every $1\leq i<j<k<l\leq n$.
    \end{enumerate}
\end{corollary}
\begin{proof}
    The implication $(i)\Rightarrow (ii)$ follows from Proposition~\ref{prop:phiro_to_slopes}.
    Every tree polynomial \eqref{eq:tree_poly}
    vanishes on $\ms_n$ by Theorem~\ref{thm:affine_slope_var}, so $(ii)\Rightarrow(iii)$.
    Now assume $(iii)$. \linebreak
    Let $\overline{\varphi}_{ij}\coloneqq\varphi(i,j)\varphi(i,j)$ be the non-normalized phases, and let $\overline{\varphi}\coloneqq (\overline{\varphi}_{ij})_{0\leq i<j\leq n}\in\C^{N}$ with $N=\binom{n+1}{2}$.
    Note that $\overline{\varphi}_{ij}\neq0$ for $i\neq j$ by uniformity of $\varphi$.
    Applying \eqref{eq:laurent_iso}, with its distinguished index $n+1$ replaced with $0$, gives
    \begin{equation}\label{eq:anchored}
        \operatorname{pen}_{ijkl0}(\overline{\varphi})
        =-(\overline{\varphi}_{0i}\overline{\varphi}_{0j}\overline{\varphi}_{0k}\overline{\varphi}_{0l})^2
        \tau_{ijkl}(\widehat{\varphi}^{\,2}) =0.
    \end{equation}
    Consequently, every pentad whose index sets contains $0$ vanishes at $\overline{\varphi}$.
    By Lemma~\ref{lem:pen_as_pen_n} all remaining pentads also vanish at $\overline{\varphi}$.
    By the five-point realizability criterion of \cite[Lemma~3.4 and Theorem~3.6]{phirotopes}, the vanishing of all these pentads implies that $\varphi$ is realizable.
    Thus $(iii)\Rightarrow (i)$.
    This finishes the proof.
\end{proof}

\printbibliography

\end{document}